\documentclass[12pt]{amsart}
\usepackage[margin=1in]{geometry}
\usepackage[utf8]{inputenc}
\usepackage{amsmath}
\usepackage{amssymb}
\usepackage{esint}%for fint

\usepackage[pdftitle={DisquesPFH4},
  pdfauthor={Samuel Bronstein},
  pdffitwindow=true,
  breaklinks=true,
  colorlinks=true,
  urlcolor=blue,
  linkcolor=red,
  citecolor=red,
  anchorcolor=red]{hyperref}
\usepackage{titleps}
\usepackage{xcolor}
\definecolor{links}{RGB}{70,0,255}
\hypersetup{citecolor=blue,linkcolor=links}

\usepackage{subfiles}

\numberwithin{equation}{section}
\numberwithin{table}{section}
\numberwithin{figure}{section}

\newcommand{\C}{\mathbb C}
\newcommand{\D}{\mathbb D}
\newcommand{\R}{\mathbb R}
\renewcommand{\H}{\mathbb H}
\renewcommand{\S}{\mathbb S}
\newcommand{\Hom}{\text{Hom}}

\newcommand{\calE}{\mathcal E}
\newcommand{\calH}{\mathcal H}
\newcommand{\calO}{\mathcal O}

\newcommand{\dbar}{\overline{\partial}}
\newcommand{\eps}{\varepsilon}
\newcommand{\Vol}{\text{Vol}}

\newcommand{\II}{\text{I\!I}}
\newcommand{\SO}{\mathrm{SO}}
\newcommand{\PSL}{\mathrm{PSL}}
\newcommand{\one}{\mathbf{1}}
\newcommand{\Isom}{\mathrm{Isom}}

\theoremstyle{plain}
\newtheorem{theorem}{Theorem}[section]
\newtheorem*{theorem*}{Theorem}
\newtheorem{lemma}[theorem]{Lemma}
\newtheorem{proposition}[theorem]{Proposition}

\newtheorem{corollary}[theorem]{Corollary}

\theoremstyle{definition}
\newtheorem{definition}[theorem]{Definition}

\newtheorem{notation}[theorem]{Notation}

\theoremstyle{remark}
\newtheorem{remark}{Remark}[section]

\title{Almost-fuchsian structures on disk bundles over a surface}
\author{samuel bronstein}
\date{2023}

\begin{document}
\begin{abstract}
	Considering an integer $d>0$,
	we show the existence of convex-cocompact representations of surface groups
	into $\SO(4,1)$	admitting an embedded minimal map with curvatures in $(-1,1)$
	and whose associated hyperbolic $4$-manifolds are disk bundles of degree $d$
	over the surface, provided the genus $g$ of the surface is large enough.
	We also show that we can realize these representations as complex variation of Hodge
	structures.
	This gives examples of quasicircles in $\S^3$ bounding superminimal disks in $\H^4$
	of arbitrarily small second fundamental form.
	Those are examples of generalized almost-Fuchsian representations which are not
	deformations of Fuchsian representations.
\end{abstract}
\maketitle
\tableofcontents
\newpage
\section{Introduction}
\subsection{Context}

\subsubsection{Minimal surfaces and convex-cocompact representations}
In this paper, we consider representations of a closed surface group into~$\mathrm{SO}(4,1)\approx\Isom(\H^4)$
and minimal surfaces in the hyperbolic $4$-space equivariant under the group action.
We prove the existence of "exotic" minimal surfaces of small curvature, that is minimal surfaces 
with arbitrarily small second fundamental form while the quotient of the hyperbolic space by the
group action is a nontrivial disc bundle over the surface considered.

Minimal surfaces are surfaces which are locally critical points of the area functional.
Considering an immersion from a Riemann surface $\Sigma$ to $\H^4$,
the image of an immersion is a minimal surface if the immersion is conformal and harmonic.
%\begin{theorem}
%Fix $\delta>0$ and $d\geq 0$. There is $g_0\geq 2$ such that, for any $g\geq g_0$,
%there is a representation $\rho:\pi_1\Sigma_g\rightarrow SO(4,1)$, and an immersion $f:\D^2\rightarrow\H^4$ satisfying:
%\begin{itemize}
%\item[(a)]
%$\rho$ is faithful, discrete, convex-cocompact.
%\item[(b)]
%$f$ is a $\rho$-equivariant embedding, minimal and superminimal.
%\item[(c)]
%The second fundamental form of $f$ satisfies $\II_f<\delta$ everywhere.
%\item[(d)]
%The hyperbolic $4$-manifold $\H^4/\rho$ is a degree $d$ disc bundle over $\Sigma_g$.
%\end{itemize}
%\end{theorem}

The study of minimal surfaces in hyperbolic manifolds is a rich and diverse topic. Schoen and Yau
\cite{SY79} and Sacks--Uhlenbeck \cite{SU82} showed that in any isotopy class of an incompressible surface in a compact manifold,
there is a weakly branched immersion whose image is minimal. See \cite{MY82} for the case when
the target manifold has mean convex boundary.
In the case of a hyperbolic three-dimensional manifold with convex boundary,
there is an immersion,
who will be an embedding if we're in the isotopy class of an incompressible surface.

Considering a faithful and discrete representation $\rho$ in $\Isom(\H^n)$,
it is then natural to ask whether the quotient manifold $\rho\backslash\H^n$ contains
a non-empty compact convex set.
When that is the case, we say that $\rho$ is \emph{convex-cocompact}.
In particular, for a convex-cocompact representation $\rho$ from a closed surface group
in $\Isom(\H^n)$, there is always an equivariant weak immersion $\D^2\rightarrow\H^4$
whose image is minimal (but has branch points in general).

While branch points can be ruled out for least area minimal surfaces in hyperbolic $3$-manifolds
with convex boundary, it is not true from the dimension 4 on.
In this regard, our theorem gives examples of minimal maps into a disc bundle over a surface
which are embeddings, so without branch points and self-intersections.

\subsubsection{Hyperbolic structures on disc bundles}
Among the family of hyperbolic $4$-manifolds whose fundamental group is a surface group,
one family is of specific interest to us: disc bundles over a Riemann surface.
Considering a Fuchsian representation of a surface group into $\SO(4,1)$,
it preserves a totally geodesic disc in $\H^4$ and the quotient manifold identifies with the
normal bundle to the quotient of the disc. Topologically the quotient manifold is the product of
a disc with the Riemann surface.
However, these are not the only examples of hyperbolic disc bundles.
Gromov--Lawson--Thurston \cite{GLT88} were the
first to give examples of hyperbolic structures on nontrivial disc bundles over a surface,
considering equivariant piecewise-linear embeddings of the disc into $\H^4$. Kuiper \cite{Kui88} 
also gave examples of those, using tessellations of the hyperbolic $4$-space.

Topologically, disc bundles over a Riemann surface are classified by the Euler characteristic
of the surface and the Euler class of the bundle. The precise list of disc bundles over a closed
surface admitting a hyperbolic structure is still an open problem. Gromov--Lawson--Thurston
conjecture
that a necessary condition to admit a hyperbolic structure would be that the degree $d$ is
smaller in absolute value than the Euler characteristic $2g-2$. Kapovich \cite{Kap89}
showed that for $0\leq d\leq g-1/11$, there is indeed a hyperbolic structure on the corresponding
disc bundle.

Our result gives a new method of finding some of those hyperbolic structures, with more geometric
properties, although we don't have an explicit condition on how small the degree has to be with
regard to $g$. Check Ho \cite{Ho14} for the latest results concerning Gromov--Lawson--Thurston's
conjecture.

\subsubsection{Almost-Fuchsian representations}

Another important question regarding minimal surfaces is the potential uniqueness in a fixed
isotopy class. Already in dimension $3$, there are examples of hyperbolic $3$-manifolds
containing arbitrarily many isotopic minimal surfaces, see Huang and Wang \cite{HW15}.
However, we may get uniqueness if we add some conditions on the manifold:
In 1983, Uhlenbeck \cite{Uhl83} introduced the notion of \emph{almost-fuchsian representation}
which are faithful and discrete representations of a surface group into $\PSL(2,\C)$ admitting
an equivariant immersion~$\D^2\rightarrow\H^3$ minimal and with principal curvatures in~$(-1,1)$.

Remark that asking the principal curvatures to be zero is equivalent to ask that our immersion
is totally geodesic.
Here as we relaxed this closed condition into asking the principal curvatures to be in  $(-1,1)$,
the notion is open in the character variety.

This differential geometric assumption has several consequences: the representation is then
convex-cocompact, the map $f$ is the only minimal map in its isotopy class and
it is an embedding.

Since then, the notion of almost-fuchsian representation has been extended to almost-fuchsian
immersions \cite{KS07}. Seppi \cite{Sep16} studied almost-fuchsian discs in $\H^3$,
Donaldson \cite{Don03}, Hodge \cite{Hod05} and
Trautwein \cite{Tra19} considered the hyperkähler structure on the moduli space of
almost-fuchsian representations extending the hyperkähler structure on Teichmüller space.
With Smith \cite{BS23}, we have exhibited a parametrization of the set of
almost-fuchsian discs in $\H^3$ by a convex set of bounded holomorphic quadratic
differentials on the disc.
There are important questions related to foliations by surfaces of almost-fuchsian manifolds,
see \cite{GHW10} and \cite{CMS23}.

Jiang \cite{Jia21} considered almost-fuchsian discs in $\H^n$. For faithful and discrete
representations into $\SO(n,1)$, the existence of an equivariant minimally immersed disc whose
scalar second fundamental forms have curvatures in $(-1,1)$ actually implies the uniqueness of
the minimal surface in its isotopy class, its embeddedness, and the convex-cocompactness of the
representation.

Recently, Davalo \cite{Dav23} extended the notion of almost-fuchsian immersion into
nearly-geodesic immersion into a higher rank symmetric space of noncompact type. The only known
examples of nearly-geodesic immersions are deformations of totally geodesic immersions.
While in $\PSL(2,\C)$ those are indeed the only examples, we bring here examples of such
immersions which cannot be deformed to totally geodesic ones.
In particular, the limit set of our representations is a quasi-circle, there is a quasiconformal
homeomorphism of $\S^3$ sending our limit set to a round curve, but this homeomorphism cannot
be made equivariant.

\subsubsection{Superminimal maps and complex variations of Hodge structures}

The construction we will exhibit will naturally give examples of complex variations of Hodge
structures in~$\SO(4,1)$.
The notion of complex variation of Hodge structures for representations of surface groups
come from the Nonabelian Hodge Correspondence,
see for instance this survey by Wentworth \cite{FMS+16}.
The Nonabelian Hodge Correspondence is a $1$ to $1$  correspondence
between the moduli space of representations of surface group with the moduli space of an object
called Higgs bundle, on which there is an action of $\C^\ast$.
The fixed points of this action are called complex variations of Hodge structures.
Loftin--McIntosh \cite{LM19} give an extensive description of those
representations in their article about equivariant minimal surfaces in $\H^4$.
For a representation $\rho$ with an equivariant minimal immersion~$f$, asking that~$\rho$
is a complex variation of Hodge structure is equivalent to asking that $f$ is superminimal,
i.e. that it has circular ellipse of curvature, see Definition 2.5.

The term superminimal was coined by Kommerell in \cite{Kom97}.
In particular, it implies the minimality of the image of the map.
Regarding the study of superminimal surfaces in the sphere $\S^4$, there are some impressive results:
see for instance Bryant \cite{Bry82}, who exhibited a duality between superminimal surfaces in $\S^4$
and holomorphic legendrian curves in $\mathbb C P^3$. Check Forstneri\v{c} \cite{For21} for the analogous
statement in $\H^4$.
%Interestingly, we don't know if the existence of a superminimal immersion can be translated
%on some property on the foliation of $\S^3$ given by the normal bundle to our immersion with its 
%limit set.

More precisely, Loftin--McIntosh show that for any pair of integers $d$ and $g\geq 2$
such that~$|d|<2g-2$,
there are representations from $\pi_1\Sigma_g$ into $\SO(4,1)$ with an equivariant superminimal
immersion whose normal bundle is of degree $d$ over $\Sigma_g$.
Our result states that for $d$ very small before $g$, there are such pairs $(\rho,f)$
where $f$ is actually embedded and $\rho$ is convex-cocompact. We also have, as the normal
bundle to $f$ will identify to the quotient manifold~$\rho\backslash\H^4$, hyperbolic structures on
degree $d$ disc bundles over $\Sigma_g$.

\subsubsection{Solving "mixed" Toda systems}%To rewrite
Finally, the method used to build our examples is to solve a non-linear PDE system
which we call a mixed Toda system.
Toda systems, and Toda systems with opposite signs are defined by Guest--Lin \cite{GL14}
and designate PDE systems where the unknown is $(w_i)$
and satisfy an equation of the type
\begin{equation}
	\Delta w_i = e^{w_i-w_{i-1}}-e^{w_i-w_{i+1}}
\end{equation}
While our system definitely cannot be written in this form, it belongs to a class
of equations arising frequently when dealing with flat curvature equations.

This kind of system appears  quite easily in the theory of flat bundles. Li--Mochizuki \cite{LM20}
encounter Toda systems (of the opposite type) when dealing with cyclic Higgs bundles associated to
representations in rank $2$ Lie groups.
Malchiodi--Ruiz \cite{MR11} and Battaglia--Jevnikar--Malchiodi--Ruiz \cite{BJM+15} made use of
the Moser--Trudinger inequality to solve some Toda systems motivated by the non-abelian
Chern-Simons theory.
Here in rank $1$ with the complex Hodge variation structure
assumption, we get something alike, but with tweaked signs:
\begin{equation}\left\{\begin{array}{lc}
	\Delta u &=-1+e^{2u}+e^{-2u}e^{2v}|\alpha|^2\\
\Delta v &=\frac{1}{2g-2}-e^{-2u}e^{2v}|\alpha|^2\end{array}\right.
\end{equation}
We think the method used in this article to solve this system could be applied to
a broader class of PDE systems, yet to be described.
\subsection{Statement of the results}
First, we state our result in dimension $4$
concerning degree $1$-disc bundles over Riemann surfaces.
\begin{theorem}
	Fix $r\in(0,1)$.
	There is $g_0\geq 2$ such that, for any $g\geq g_0$
	and $\Sigma$ a closed Riemann surface of genus $g$,
	there is a representation $\rho:\pi_1\Sigma\rightarrow\Isom(\H^4)$
	and an immersion $f:\S^2\rightarrow\H^4$ satisfying:
	\begin{itemize}
		\item[(a)]
		$\rho$ is faithful, discrete, convex-cocompact, almost-fuchsian
		\item[(b)]
		$f$ is a $\rho$-equivariant embedding, superminimal
		\item[(c)]
		The second fundamental form of $f$ satisfies $|\II_f|<r$ everywhere.
		\item[(d)]
		The quotient hyperbolic $4$-manifold $\rho\backslash\H^4$
		is a degree $1$ disc bundle over $\Sigma$.
	\end{itemize}
\end{theorem}
As a corollary, we have the existence of almost-fuchsian representations uniformizing 
degree $d$ disc bundles over Riemann surfaces of high genus.
\begin{corollary}
	Let $d>0$. For any $g_0>0$, there is $g\geq g_0$, $\Sigma$ a surface of genus $g$
	; $\rho:\pi_1\Sigma\rightarrow\Isom(\H^4)$ and $f:\D^2\rightarrow\H^4$satisfying:
	\begin{itemize}
	\item[(a)]
	$\rho$ is faithful, discrete, convex-cocompact, almost-fuchsian
	\item[(b)]
	$f$ is a $\rho$-equivariant embedding, superminimal.
	\item[(c)]
	The second fundamental form of $f$ satisfies $\II_f<r$ everywhere.
	\item[(d)]
		The quotient manifold $\rho\backslash\H^4$ is a degree $d$ disc bundle over $\Sigma$.
	\end{itemize}
\end{corollary}

\subsection{Scheme of proof}
\subsubsection{Writing the Flat curvature equations}
In order to build a hyperbolic structure on a degree $d$ disc bundle, we first show
the existence of a hyperbolic structure on a degree $1$ disc bundle over a surface with the
desired properties. Once this example is built, we can take a cover of degree $d$ of this bundle
in order to obtain a hyperbolic structure on a degree $d$ disc bundle over a surface (of genus
large enough).

Consider a Riemann surface $\Sigma$ of genus $g$.
Our proof will rely on solving a PDE system of the following kind:
\begin{equation}
	\left\{\begin{array}{cl}
		\Delta u &= -1 +e^{2u}+e^{-2u}e^{2v}|\alpha|^2\\
		\Delta v &= \frac{1}{2g-2}-e^{2v}e^{-2u}|\alpha|^2
	\end{array}\right.
\end{equation}
where $\alpha$ has to be a holomorphic section of a degree $4g-3$ line bundle over $\Sigma$.
The Laplacian is always taken with respect to the Poincaré volume form on the surface,
and every line bundle will be considered endowed with a metric whose curvature form
is proportional to the Poincaré volume form, hence defining the norm $|\alpha|$.

The second section will be devoted to showing how solving this system yields a minimal disc in $\H^4$ equivariant under a representation $\rho$. With the good pointwise control
on~$e^{2v}e^{-2u}|\alpha|^2$, the representation $\rho$ will be almost-fuchsian
and the quotient manifold $\rho\backslash\H^4$ will be a degree $1$ disc bundle over the Riemann surface.

The associated minimal map will have a holomorphic second fundamental form (see Definition 2.2.)
$\alpha\oplus 0$, hence be superminimal.
\subsubsection{Controlling the geometric data}
The third section will be devoted to finding families~$(\Sigma_n)$ of surfaces with genus going to infinity,
with data $(\alpha_n)$ associated,
while controlling the behavior of solutions of PDEs on $(\Sigma_n)$ depending on $(\alpha_n)$.
To do so, we need to control the systole and spectral gap of $(\Sigma_n)$,
and the ratio $\frac{\fint|\alpha_n|^2}{|\alpha_n|_\infty^2}$.

Controlling the systole is the easiest part, as the systole is nondecreasing when taking
the Riemannian cover of a surface.
Controlling the spectral gap is trickier, but can be done by taking $(\Sigma_n)$ a random cover
of a genus $2$ surface thanks to a result of
Magee--Naud--Petri \cite{MNP22}:
considering a cover $\Sigma'\rightarrow\Sigma$ of closed Riemann surfaces,
the relative spectrum of $\Sigma'$ is defined to be the part of the spectrum of~$\Sigma'$
which is not in the spectrum of $\Sigma$.
Magee--Naud--Petri states that for a random cover of large degree of $\Sigma$,
the probability that the relative spectral gap is lower than $\frac{3}{16}-\eps$ goes to zero
when the degree goes to infinity.
Using this result, we are able to find families
of coverings of a Riemann surface whose spectral gap is uniformly lower bounded.

Finally, we need to build $\alpha_n$ such that the ratios
$\frac{|\alpha_n|_\infty^2}{\fint|\alpha_n|^2}$ are uniformly bounded.
This property will be called \emph{balanced} family of sections.
We will show that multiplying a cover of a section with a section of a line bundle of degree $1$
yields a balanced family of sections of the desired type.

\subsubsection{Building the maps $\Phi$ and $\Psi$}
Once this is done, we need to solve our PDE system.
For this, we will consider each line of the system and build two maps.
Fixing $v$ and solving the first line in $u$, we will build a map
$\Phi:C^{0,\alpha}(\Sigma)\rightarrow C^{2,\alpha}(\Sigma)$ which will be well-defined on the
set $\mathcal K_v$, where:
\begin{equation}
	\mathcal{K}_v:=\{v\in C^{0,\alpha}(\Sigma): e^{2v}|\alpha|^2\leq\frac{1}{2}\}
\end{equation}
And it happens that the image of $\mathcal{K}_v$ is in the set $\mathcal{K}_u$, where
\begin{equation}
	\mathcal{K}_u:=\{u\in C^{0,\alpha}(\Sigma): -\frac{\ln 2}{2}\leq u\leq 0\}
\end{equation}
We would like to do the same work for the second equation, fixing $v$, solving in $u$ to get
a map $\Psi$. However, the work to build $\Psi$ is considerably harder.
The third section will be devoted to the study of $\Psi$. The method used to solve it is heavily
inspired from the analysis of the prescribed curvature equation in the projective plane
$\R P^2$: for an even function on the sphere which is nonnegative on an open set, it is
the curvature of a metric conformal to the round metric on $\S^2$ (and even). See
Chang--Yang \cite{CY03} for a precise explanation of this.
With the same tools, we are able to solve the second line and to get the following estimate:
There is a constant $C>0$
such that the solution $v$ satisfies:
\begin{equation}
	e^{2v}e^{-2u}|\alpha|^2\leq\frac{C}{2g-2}
\end{equation}
The involved constant $C$ will depend on three parameters:
$\delta$ the systole of the surface, $\Lambda$ its spectral gap and
the ratio between the $L^\infty$-norm of $\alpha$ and its $L^2$-average.
Hence thanks to the third section it can be bounded on surfaces of arbitrarily large genus.

\subsubsection{Piecing everything together}
Considering a family of covers $(\Sigma_n)$ with a balanced family of sections of sections $(\alpha_n)$
and controlled systole and spectral gap, we are able, provided~$n$ is large enough,
to build a map~$\Psi:\mathcal{K}_u\rightarrow\mathcal{K}_v$,
which will be continuous.

With $\Psi$ constructed and well-defined,
We will show that the composed map $\Phi\circ\Psi$ defined on $\mathcal K_u$
has image in a compact convex subset of $\mathcal K_u$.
Hence by Schauder's fixed-point theorem,~$C^{0,\alpha}(\Sigma)$ being a Banach space,
$\Phi\circ\Psi$ has a fixed point $u\in\mathcal K_u$.
Denote $v=\Psi(u)$.

It is clear that $(u,v)$ is a solution of the PDE system considered,
and so it corresponds to an almost-fuchsian representation $\rho$
whose quotient manifold is a degree $1$ disc bundle over the Riemann surface.

The author thanks Nicolas Tholozan (DMA) for his help and guidance during
this work.

\section{The curvature equations}
\subsection{Minimal equivariant immersions in $\H^4$}
Let $\Sigma$ be a closed oriented surface of genus~$g\geq 2$.

Consider a representation $\rho:\pi_1\Sigma\rightarrow\SO_0(4,1)\approx\Isom_+(\H^4)$.
We will always assume that~$\rho$ is faithful and discrete.
Let $f:\D^2=\tilde\Sigma\rightarrow\H^4$ be an equivariant immersion with minimal image.

Endowing $\Sigma$ with the induced metric by $f$,
we get a conformal and harmonic map.

Associated to $f$, there is a flat $\R^{4,1}$ vector bundle $E$,
which decomposes as follows:

Denote $\II$ the second fundamental form of $f$,
and $B$ its shape operator.
The vector bundle~$(E,\nabla_E)$ splits as $E=O_\R\oplus T_\R \oplus N_\R$,
where $O_\R$ is the trivial vector bundle of rank $1$,
$T_\R$ the tangent bundle is a rank $2$ vector bundle,
and $N_\R$ is the normal bundle to $f$, again of rank $2$.
In this splitting, the connection has the following expression:
\begin{equation}
\nabla_E=\left(\begin{array}{ccc}
\nabla & a^\dagger & 0\\
a & \nabla_T & B \\
0 & \II & \nabla_N \end{array}\right)
\end{equation}
Denote $q$ the quadratic form in $\R^{4,1}$ given by
\begin{equation}
q(x_1,\ldots,x_5):=x_1^2+x_2^2+x_3^2+x_4^2-x_5^2
\end{equation}
As the map $f$ is valued in the quadric $q=-1$ in $\R^{4,1}$,
the image of $f$ gives a reduction of~$E$ into a $S(O(1)\times O(4))$ bundle.
The differential of $f$ being valued, by definition, in the Tangent bundle,
we denote $a=df$ the nonvanishing section of $T^\ast\Sigma\otimes\Hom(O_\R,T)$.
The definition of $\II$ and of $B$ gives the desired decomposition of $\nabla_E$.

\begin{remark}
This splitting is orthogonal with respect to $b$ the signature $(4,1)$ bilinear form,
with $O_\R$ being a negative line.
As the surface is oriented and our representation is into orientation-preserving
isometries $\SO_0(4,1)\approx\Isom_+(\H^4)$,
$T_\R$ and $N_\R$ are endowed with an~$SO(2)$ structure, corresponding to the
induced metrics on the tangent and normal bundle to $f$.
The nondegenerate metric $b$ also implies that $\II$ and $B$ are antiadjoint with respect to~$b$:
\begin{equation}
	b(B(w)u,v)+b(\II(u,v),w)=0\quad\forall (u,v)\in T_\R,\, w\in N_\R\,. 
\end{equation}
We denote $B=-\II^\dagger$,
where $\dagger$ denotes the transposition,
identifying $\Hom(E,F)$ with $\Hom(F^\ast,E^\ast)$.
\end{remark}
Complexify now the bundles $T_\R$ and $N_\R$.
As each of these are equipped with an $\SO(2)$-structure,
their complexification splits into the sum of two line bundles,
the restriction of~$B$ giving a nondegenerate pairing between those.
Now for each line subbundle of $E$ coming into play,
we consider $\nabla_L$ the induced connection of $\nabla_E$
and we take its $(0,1)$-part to get~$\dbar_L$ a holomorphic structure on $L$.
Remark that the $(0,1)$-part is well-defined because the induced metric
endows the surface with a Riemann surface structure.

The diagonalization of the $SO(2)$ structure with regard to the metrics yields:
\begin{proposition}
	Viewing $(\II_\C)$ as a $1$-form valued in $\Hom(K^{-1},N)$,
the complexification of the second fundamental form can be written:
\begin{flalign}
(T_\R)_\C&=\overline K\oplus K=K^{-1}\oplus K\\
(N_\R)_\C&=N\oplus\overline N=N\oplus N^{-1}\\
(\II)_\C&=\left(\begin{array}{cc}
\alpha & \beta^\ast \\ \beta & \alpha^\ast \end{array}\right)
\end{flalign}
with $\alpha$ and $\beta$ $(1,0)$-forms.
\end{proposition}
\begin{proof}
The splitting of the complexification of an $\SO(2)$-bundle is given by the eigenspaces of $J$
the fiberwise rotation of angle $\frac{\pi}{2}$.
Remark that this splitting diagonalizes $\nabla_T$ and $\nabla_N$.
In this splitting, write down the second fundamental form:
\begin{equation}
(\II)_\C=\left(\begin{array}{cc}
\alpha &\gamma \\ \beta & \delta \end{array}\right)
\end{equation}
As $\II$ is real-valued, we have:
\begin{flalign}
	\gamma&=\overline \beta=\beta^\ast\\
	\delta&=\overline \alpha=\alpha^\ast
\end{flalign}
Finally, as $f$ is minimal, $\II$ is traceless, and so the $(1,1)$-part of the second fundamental
form is zero:
This implies that  both $\alpha\in \Gamma(K\otimes\Hom(K^{-1},N))$ and $\beta\in\Gamma(K\otimes\Hom(K^{-1},N^{-1}))$ are $(1,0)$-forms.
\end{proof}
\begin{remark}
	In the process, we endowed the bundle $E\otimes\C$ with the following holomorphic
	structure:
	\begin{flalign}
		E\otimes\C&=\calO\oplus K^{-1}\oplus K\oplus N\oplus N^{-1}\\
		\dbar_E&=\dbar_\calO\oplus\dbar_{K^{-1}}\oplus\dbar_K\oplus\dbar_N\oplus\dbar_{N^{-1}}
	\end{flalign}
\end{remark}
It is clear that the sections $\alpha,\beta$ characterize wholly the second fundamental form.
Hence we define the \emph{Holomorphic second fundamental form}:
%We recall our definition of the Hopf differential:
\begin{definition}
	Let $f:\D^2\hookrightarrow\H^4$ be a minimal immersion.
	The \emph{Holomorphic second fundamental form} of $f$
	is the $(2,0)$-part of the complexification of the second fundamental form:
	\begin{equation}
		\II^{2,0}=\alpha\oplus\beta\in \Gamma(K^2(N\oplus N^{-1}))
	\end{equation}
\end{definition}

\subsection{Writing the flat curvature equations as a scalar PDE system}
In order to write this system as a scalar PDE system, we need to write down each metric as
a conformal change of a preferred one.
For the metric on $K$, we denote $e^{2u}\sigma$ the induced metric, where~$\sigma$
is the unique complete hyperbolic metric in the conformal class.

Denote $\omega$ the Poincaré $(1,1)$-form associated to $\sigma$, which is locally $dz\wedge d\overline z$.
Also, we denote by $\Delta$ the Laplace-Beltrami operator on $(\Sigma,\sigma)$.
Note that $\Delta$ is defined to be locally $\partial_x^2+\partial_y^2$, so it is a nonpositive
operator satisfying:
\begin{equation}
	\forall (f,g)\in C(\Sigma),\quad\int_\Sigma(\Delta f)g \,\omega
	+\int_\Sigma \nabla f\cdot\nabla g\,\omega=0
\end{equation}

Generally, on a holomorphic line bundle $L$ of degree $d$,
there is a hermitian metric $h_L$ whose curvature form
is the following:
\begin{equation}
	\frac{i}{2\pi}F_{L,h_L}=\frac{d}{2g-2}\omega
\end{equation}
The metric $h_L$ is uniquely defined up to multiplicative constant.
From now on, every line bundle will be endowed with a metric of this kind,
and a given fiberwise norm $|\cdot|_L$.

Let's write down the flat curvature condition as a scalar PDE system:
\begin{theorem}
Let $f:\Sigma\hookrightarrow M$ be a conformal harmonic immersion into a complete hyperbolic
$4$-manifold.
Consider $\alpha\oplus\beta$ its holomorphic second fundamental form.
Denote $e^{2u}\sigma$ the induced metric on $\Sigma$, $e^{2v}h_N$ the induced metric on the normal
bundle.
Then the flat curvature equations are equivalent to the system:
\begin{equation}
	\left\{\begin{array}{cl}
\dbar \alpha &= 0\\
\dbar \beta &=0\\
\Delta u &=e^{2u}-1+ e^{-2u}(e^{2v}|\alpha|^2+e^{-2v}|\beta|^2\\
\Delta v &=c-e^{-2u}(e^{2v}|\alpha|^2-e^{-2v}|\beta|^2)
	\end{array}\right.
\end{equation}
\end{theorem}
\begin{proof}
First, as $B=-\II^\dagger$, 
we have the explicit description of $B$
\begin{equation}
B=\left(\begin{array}{cc}
-\alpha^{\dagger,\ast} &-\beta^{\dagger,\ast}\\
-\beta^\dagger & -\alpha^\dagger\end{array}\right)
\end{equation}
Now the flat curvature equation is $F_{\nabla_E}=0$, so coordinatewise:
\begin{equation}
\left\{\begin{array}{rl}
\nabla \alpha&=0\\
\nabla \beta&=0\\
	F_{K^{-1}}  + a a^\ast-\alpha^{\ast,\dagger}\alpha-\beta^{\ast,\dagger}\beta&=0\\
	F_N -\alpha\alpha^{\ast,\dagger}-\beta^\ast\beta^\dagger&=0
\end{array}\right.
\end{equation}
where we wrote $\alpha\alpha^\ast$ for the wedge product $\alpha\wedge\alpha^\ast$.
Also, the notation $\nabla\alpha=0$ means that~$\alpha$ is parallel for the connection on $\Hom(K^{-1},N)$, namely $\nabla_N\wedge\alpha +\alpha\wedge\nabla_{K^{-1}}=0$.

The standard formula for the curvature under a conformal change of metric is
\begin{equation}
	\left\{\begin{array}{ll}
iF_{K^{-1}}&= e^{-2u}(-\Delta u-1)\omega\\
iF_N&= e^{-2u}(-\Delta v +c)\omega
\end{array}\right.\end{equation}
Now recall that as $\alpha$ is a $(1,0)$-form, and locally $dz\wedge d\overline z=-2idx\wedge dy$,
hence 
\begin{equation}
\left\{\begin{array}{cl}
iaa^\ast&=\omega\\
i\alpha\alpha^\ast&=e^{-4u}e^{2v}|\alpha|^2 \omega\\
i\beta\beta^\ast &=e^{-4u}e^{-2v}|\beta|^2 \omega
\end{array}\right.
\end{equation}
So the flat curvature equation is equivalent to:
\begin{equation}
	\left\{\begin{array}{cl}
e^{-2u}(\Delta u +1)&=1+e^{-4u}(e^{2v}|\alpha|^2+e^{-2v}|\beta|^2)\\
e^{-2u}(\Delta v -c)&=e^{-4u}(-e^{2v}|\alpha|^2+e^{-2v}|\beta|^2)
	\end{array}\right.\end{equation}
Also, as $\alpha$ is of type $(1,0)$, necessarily $\partial\alpha=0$.
Hence $\nabla \alpha=0$ rewrites as $\dbar\alpha =0$, and the same applies to $\beta$.
Finally, rearranging the right and left hand sides, we get the desired scalar PDE system:
\begin{equation}
	\left\{\begin{array}{cl}
\dbar\alpha &=0\\
\dbar\beta &=0\\
\Delta u &= -1+e^{2u}+e^{-2u}(e^{2v}|\alpha|^2+e^{-2v}|\beta|^2)\\
\Delta v &= c-e^{-2u}(e^{2v}|\alpha|^2-e^{-2v}|\beta|^2)
	\end{array}\right.
\end{equation}
\end{proof}

\begin{remark}
	The Higgs bundle associated to $f$ can be read with our notation, it will be the
	following: $\calE=\calO\oplus K^{-1}\oplus K\oplus N\oplus N^{-1}$,
	and the holomorphic structure with the Higgs field have the explicit description:
	\begin{equation}
		\dbar_\calE=\left(\begin{array}{ccccc}
			\dbar & 0&0&0&0\\
			0&\dbar&0&-\alpha^\ast&-\beta^\ast\\
			0&0&\dbar&0&0\\0&0&\beta^\ast&\dbar&0\\
		0&0&\alpha^\ast&0&\dbar\end{array}\right),\quad\Phi_\calE=
		\left(\begin{array}{ccccc}
			0&a&0&0&0\\0&0&0&0&0\\a^\ast&0&0&0&0\\0&0&0&0&0\\0&0&0&0&0
		\end{array}\right)
	\end{equation}
	Recall that we chose $N$ so that $\deg N\geq 0$.
	Observing that $\calO\oplus K^{-1}\oplus N$ is an invariant subbundle,
	the stability condition implies:
	\begin{equation}
		0\leq\deg N\leq 2g-2\,.
	\end{equation}
\end{remark}
\begin{proposition}
	Let $f$ be a minimal map with $\alpha=0$.
	Then $\deg N=0$ and $f$ is totally geodesic.
\end{proposition}
\begin{proof}
	As $\alpha=0$, $N$ is an invariant subbundle of the corresponding Higgs bundle to $f$.
	Then by semi-stability, $\deg N=0$.
	Now if $\beta\neq 0$, $v$ is solution of:
	\begin{equation}
		\Delta v =e^{-2u}e^{-2v}|\beta|^2
	\end{equation}
	Integrating this identity implies
	\begin{equation}
		\int_\Sigma e^{-2u}e^{-2v}|\beta|^2=0
	\end{equation}
	Hence $\beta$ is the zero section and $f$ is totally geodesic, as asserted.
\end{proof}
We will be interested in the specific case when the second fundamental form,
viewed as a map $T_\R\rightarrow (T_\R)^\ast$, has \emph{circular ellipse of curvature},
i.e. the image of the unit circle in $T_\R$ by the map $X\mapsto\II(X,X)$ is a circle.
This corresponds to asking the $(2,0)$ and $(0,2)$ part of $\II_\C$ to be orthogonal.

Following Loftin--McIntosh \cite{LM19}:
\begin{definition}
We say that $f$ is superminimal if $f$ is harmonic, conformal and
its Holomorphic second fundamental form $\alpha\oplus\beta$ satisfies:
\begin{equation}
\alpha\cdot\beta=0\in H^0(K^4)
\end{equation}
\end{definition}
\begin{remark}
	The notion of superminimal maps goes back to 1897 with Kommerell \cite{Kom97} who studied
	superminimal immersions in $\S^4$. See Forstneri\v{c} \cite{For21} for a recent review
	of the topic.
	Another interpretation of the superminimality is that it is
	equivalent to 
	asking the normal Gauss map, valued in the Grassmannian of
	geodesic disks in $\H^4$, is conformal.
\end{remark}
Depending on the degree of $N$ superminimal maps behave differently:
\begin{proposition}
Let $f$ be a superminimal map with $\deg N=0$,
equivariant under a faithful and discrete representation $\rho:\pi_1\Sigma\rightarrow SO(4,1)$,
$\Sigma$ being a closed surface.
Then $f$ is totally geodesic and $\rho$ is a fuchsian representation.
\end{proposition}
\begin{proof}
While this proposition can be found in \cite{LM19},
we reproduce it here as it can be seen as a consequence of our PDE system:
As the normal bundle is trivial, we have $c=0$ and the superminimality implies $\alpha=0$ or
$\beta=0$. For convenience's sake, assume $\beta=0$.
Then~$(u,v)$ must be a solution of the system
\begin{equation}\left\{\begin{array}{cl}
\Delta u &= -1+e^{2u}+e^{-2u}e^{2v}|\alpha|^2\\
\Delta v &= -e^{2u}e^{2v}|\alpha|^2\end{array}\right.
\end{equation}
Fix $u$, and observe that, on a closed Riemann surface, we must have $\int \Delta v =0$.
Hence  $\alpha=0$, the immersion $f$ is totally geodesic and $\rho$ is fuchsian.
\end{proof}
\subsection{Almost-fuchsian disks in $\H^4$}
First introduced by Uhlenbeck for immersions into~$\H^3$ \cite{Uhl83},
the notion of almost-fuchsian immersion naturally generalizes to immersions into~$\H^4$.
\begin{definition}
	Let $f:\D\hookrightarrow \H^4$ be a proper immersion.
	$f$ is said to be \emph{almost-fuchsian} if it is minimal and
	there is $\delta>0$ such that the second fundamental form
	satisifes:
	\begin{equation}
		\forall u\in T_f,\,|\II(u,u)|\leq (1-\delta)|u|^2
	\end{equation}
\end{definition}
\begin{proposition}
	Let $f:\D\rightarrow \H^4$ be a proper almost-fuchsian immersion.
	Then the exponential map $N_f\rightarrow\H^4$ is a diffeomorphism
\end{proposition}
As this constitutes, with the following result, the lemma 2.2. of \cite{Jia21}, we won't reproduce
a proof here. Note that Jiang's convention for the second fundamental form differs from ours by a
factor one half, hence the difference in the condition.
%\begin{proof}
%	Consider the exponential map $\exp_f:\D^2\times\R^2\rightarrow\H^4$.
%	As $f$ is proper, the exponential map is surjective.
%	
%	We will write down the induced metric by $\exp_f$, and show that it is nondegenerate.
%	In polar coordinates, as the radial direction is geodesic,
%	Denote $B_f=(\II_f)^\perp$ the shape operator of $f$.
%	we can decompose the induced metric as $\exp_f^\ast g_{\H^4}=(g_r)_{r\geq 0}$
%	with
%	\begin{flalign}
%		g_0=f^\ast g_{\H^4}
%	\end{flalign}
%	As $f$ is an immersion, $g_0$ is nondegenerate, and the extension $\tilde g_0=g_0\times g_{S^1}$
%	is a nondegenerate metric on $\D^2\times\S^1$.
%
%	Hence we define $B_t$ the endomorphism such that:
%	\begin{equation}
%		\forall t\geq 0,\,g_t(X,Y)=\tilde g_0(B_tX, B_tY)
%	\end{equation}
%	A computation in local coordinates gives the following:
%	\begin{flalign}
%		B_0&=\text{id}\\
%		\dot B_0&=B_f\times\text{id}\\
%		^t(\ddot B_t-B_t)B_t&=0
%	\end{flalign}
%	As $B_0$ and $\dot B_0$ commute, we deduce
%	\begin{equation}
%		B_t=\cosh(t)B_0+\sinh(t)\dot B_0
%	\end{equation}
%	Now the almost-fuchsian condition states that the eigenvalues of $B_f$ are in $(-1,1)$.
%	Hence one can see that $B_t$ is invertible for all $t>0$,
%	and the pullback of the metric by the exponential map is nondegenerate.
%
%	This implies that $\exp_f$ is a local diffeomorphim, and it is surjective,
%	thus it is a global diffeomorphism, as asserted.
%\end{proof}
More specifically, almost-fuchsian disks in $\H^4$ are embeddings,
bound a quasi-circle at infinity, and if they are $\rho$-equivariant
for some $\rho:\Gamma\rightarrow\text{SO}(4,1)$ discrete and faithful,
$\rho$ is convex-cocompact.
In particular,
\begin{proposition}
	Let $\Sigma$ be a closed hyperbolic surface,
	and $\rho:\pi_1\Sigma\rightarrow\text{SO}(4,1)$
	a faithful and discrete representation such that there exists
	$f:\D\hookrightarrow\H^4$ a $\rho$-equivariant proper almost-fuchsian immersion.
	Then $\rho$ is convex-cocompact,
	$f$ is an embedding and the hyperbolic $4$-manifold~$M=\rho\backslash\H^4$
	is diffeomorphic via $\exp_f$ to a disk bundle over $\Sigma$.

	Also $f$ is the unique $\rho$-equivariant minimal immersion.
\end{proposition}
\begin{definition}
	A representation $\rho:\Gamma\rightarrow\text{SO}(4,1)$.
	is said to be \emph{almost-fuchsian} if it is discrete, faithful and 
	it admits a proper almost-fuchsian $\rho$-equivariant
	immersion $f:\widetilde\Sigma\rightarrow\H^4$.
\end{definition}
It is not known which disk bundles can be uniformized by almost-fuchsian representations.
Fuchsian representations uniformize the trivial bundle,
and we will show in this paper that for a surface of large enough genus, there are almost-fuchsian
representations uniformizing vector bundles with positive degree over the surface.
\begin{proposition}\label{fundathm}
	Let $\Sigma$ be a closed hyperbolic surface of genus $g$.
	Denote $w$ its volume form.
	Consider $N$ a line bundle of degree $d\geq 0$,
	endowed with a metric $h$ of curvature form~$c\omega=\frac{d}{2g-2}\omega$.
	Let $\alpha\in H^0(K^2 N)$ be such that there exist smooth functions $u,v$ on $\Sigma$
	satisfying:
	\begin{flalign}
		\Delta u &= -1 + e^{2u}+e^{-2u}e^{2v}|\alpha|^2\\
		\Delta v &= c - e^{-2u}e^{2v}|\alpha|^2\\
		\sup e^{-4u}e^{2v}|\alpha|^2&<1
	\end{flalign}
	Then there is a convex-cocompact representation $\rho:\pi_1\Sigma\rightarrow SO(4,1)$
	and a minimal, superminimal, $\rho$-equivariant and almost-fuchsian immersion
	$f:\D\rightarrow\H^4$ such that:
	\begin{enumerate}
		\item
			The induced metric by $f$ is $e^{2u}w$
		\item
			The hyperbolic manifold $\rho\backslash\H^4$ identifies with the normal
			bundle to $f$ and has degree~$d$
		\item
			The metric induced on the normal bundle is $e^{2v}h$
		\item
			The Holomorphic second fundamental form of $f$ is $\alpha$.
	\end{enumerate}
\end{proposition}
\begin{proof}
	With the given data, consider the bundle $\calE=\calO\oplus K^{-1}\oplus K\oplus N\oplus N^{-1}$, with the metric gotten by product of the metrics on each line.
	Adding $\alpha$ and its $h$-dual $\alpha^\ast$,
	the 2 equations are exactly the flat curvature equation,
	corresponding to a representation $\rho$ and a map $f:\D\rightarrow\H^4$ whose 
Holomorphic second fundamental form  satisfies:
	\begin{equation}
		(\II_f)_\C=\alpha\oplus 0 \in \Gamma(K^2(N\oplus N^{-1}))
	\end{equation}
	It is clear that such an $f$ is minimal and superminimal.
The last hypothesis
\begin{equation}
\sup e^{-4u}e^{2v}|\alpha|^2<1
\end{equation}
ensures that $f$ and $\rho$ are almost-fuchsian, as asserted.
\end{proof}

\section{Balanced holomorphic sections of positive line bundles}
In this section, we show that there are families of holomorphic sections of line bundles
of degree $4g-3$ which behave roughly like covers in the following sense:
\begin{definition}[Balanced family]
Let $(\Sigma_g)$ be a family of genus $g$ Riemann surfaces, equipped with their hyperbolic metric.
Let $N_g$ be a family of line bundles of degree $n_g\geq 0$
Let $(\alpha_g)$ be a family of holomorphic sections of $N_g$.

On each $\Sigma_g$, denote by $\omega_g$ the Poincaré volume form.
On each $N_g$, denote $h_g$ the metric on $N_g$ with curvature form $\frac{n_g}{2g-2}\omega_g$.
We say that $(\alpha_g,N_g)$ is \emph{balanced} if there is $C>0$ such that:
\begin{equation}
\sup_g\frac{\sup_{x\in X_g}|\alpha_g|^2_{\omega_g,h_g}}{\fint_{X_g}|\alpha_g|^2_{h_g}d\omega_g}
\leq C
\end{equation}
\end{definition}
The main idea behing this control is that when solving elliptic PDEs depending
on $\alpha_g$, as $\alpha_g$ is coarsely periodic, then the solutions should
be too, and we will show that the solutions satisfy estimates which are
independent of the genus $g$.

We will make extensive use of the following notations:
\begin{notation}
In this section, $\Sigma$ denotes a closed hyperbolic surface of genus $g$.
We denote respectively by $\delta$ and $\Lambda$ the systole and spectral gap of $\Sigma$.
Any $N$ positive line bundle over~$\Sigma$ is endowed with a metric of curvature form
proportional to the Poincaré volume form $\omega$.
\end{notation}
\begin{remark}
It is easy to exhibit examples of balanced sections  for $n_g$ proportional to $g-1$:
Indeed, if $n_g=a(g-1)$, consider $N_2$ a degree $a$ line bundle over $\Sigma_2$,
and $\alpha$ a holomorphic section of $N_2$.
and consider covers $\Sigma_g\rightarrow\Sigma_2$. Then the family of lifts of $N_2$ and $\alpha$
yields a balanced family.
\end{remark}
Here we will prove the existence of balanced families for line bundles of degree $4g-3$:
\begin{theorem}\label{equisec}
Let $\Sigma$ be a closed hyperbolic surface of genus $g\geq 2$
Consider $(\Sigma_n)$ a family of covers of $\Sigma$ of degree $n$.
Then there is a balanced family $(\alpha_{n},N_{n})$,
with $N_n$ a line bundle of degree $(4g-4)n+1$ over $\Sigma_n$.
\end{theorem}
The main idea will be to consider covers to get a balanced section of a degree $4n(g-1)$
line bundle, and then to tensorize with a well-chosen degree $1$ line bundle to get the 
desired family of sections.
This requires a description of the oscillatory behavior of a degree $1$ bundle section,
which we provide here.

Before all, we will make use of the continuous embedding $W^{2,2}\hookrightarrow L^\infty$.
\begin{lemma}\label{22bdd}
	Let $\Sigma$ be a hyperbolic surface with systole $\delta$.
	Denote $M$ the constant of continuity of the embedding $W^{2,2}(\H^2)\hookrightarrow L^\infty(\H^2)$.
	Then, for any $u\in W^{2,2}(\Sigma)$
	\begin{equation}
	|u|_\infty\leq \frac{C}{\delta}|u|_{2,2}\,.
	\end{equation}
\end{lemma}
\begin{proof}
	As smooth maps are dense in $W^{2,2}(\Sigma)$, consider $u\in W^{2,2}(\Sigma)$ smooth.
	Let $z_0\in\Sigma$.
	By definition of the systole, there is an isometric embedding from a hyperbolic disk of
	radius~$\delta$ around $z_0$.
	Fix $\chi$ a smooth radial map on $D(z_0,\delta)$,$\frac{1}{\delta}$ Lipschitz,
	between $0$ and~$1$,
	such that~$\chi(z_0)=1$ and $\chi(\partial D(z_0,\delta))=0$.
	Choose for $\chi$ a rescaling of a cutoff map on~$D(0,1)\subset\H^2$, so that
	\begin{equation}
		\sup{ D\chi}=\sup{D^2\chi}=\frac{1}{\delta}
	\end{equation}

	Then we can extend by zero to consider $u\chi\in W^{2,2}(\H^2)$.
	Using our cutoff function, we get
	\begin{equation}
		|u(z_0)|=|u\chi(z_0)|\leq |u\chi|_\infty\leq M|u\chi|_{2,2}\leq \frac{M}{\delta}|u|_{2,2}
	\end{equation}
	, as asserted.
\end{proof}
\begin{proposition}\label{gbdd}
	%Let $\Sigma$ be a hyperbolic surface with systole $\delta$ and spectral gap $\Lambda$.
	%Let $N$ be a degree $1$ line bundle over $\Sigma$ and $s\in H^0(N)$ section, with a
	%unique zero $z_0$.
	Let $z_0\in\Sigma$
	Denote $D=D(z_0,\delta)$ the disk of radius $\delta$ around $z_0$.
	Denote  by~$a$ the constant $a=\frac{2\pi}{\Vol(D)}$.
	Consider $g:\Sigma\rightarrow\R$ the $C^1$ function satisfying:
	\begin{equation}
		\left\{\begin{array}{cl}
		\int_\Sigma g &=0\\
		\Delta g &= \frac{1}{2g-2}-a\one_D\end{array}\right.
	\end{equation}
	Then there is a bound $C(\delta,\Lambda)$ such that:
	\begin{equation}
		|\sup g -\inf g|\leq C(\delta,\Lambda)
	\end{equation}
\end{proposition}
\begin{proof}
	as $a=\frac{2\pi}{\Vol D}$,
	\begin{equation}
		\int_\Sigma \frac{1}{2g-2}-a\one_D=0
	\end{equation}
	hence $g$ is well-defined.
	
	First, we have $\int g\Delta g =-|\nabla g|_2^2$, so
	\begin{flalign}
		|\nabla g|_2^2=a\int_D g\leq 2\pi\sup g
	\end{flalign}
	Recall Bochner's identity
	\begin{equation}
		\int |\nabla^2g|^2=\int (\Delta g)^2+\int |\nabla g|^2
	\end{equation}
	which translates to give the estimate
	\begin{equation}
		|\nabla^2 g|_2^2\leq
		\frac{4\pi^2}{\Vol (D)}-\frac{2\pi}{2g-2}+2\pi\sup g
	\end{equation}
	As $\int g=0$, by definition of the spectral gap
	\begin{equation}
		|g|_2^2\leq\Lambda|\nabla g|_2^2\leq 2\pi\Lambda\sup g\,.
	\end{equation}
	Combining all this, we bound the $W^{2,2}$ norm
	\begin{equation}
		|g|_{2,2}^2=|\nabla^2 g|_2^2+|\nabla g|_2^2+|g|_2^2\leq
		\frac{4\pi^2}{\Vol(D)}-\frac{2\pi}{2g-2}+2\pi(2+\Lambda)\sup g
	\end{equation}
	By lemma \ref{22bdd}, this implies
	\begin{equation}\label{eqsupinf}
		|\sup g -\inf g|\leq\frac{2M}{\delta}\bigg(
		\frac{4\pi^2}{\Vol(D)}-\frac{2\pi}{2g-2}+2\pi(2+\Lambda)\sup g\bigg)^{\frac{1}{2}}\,.
	\end{equation}
	Because $\int g=0$, we know that $\inf g\leq 0$ and $\sup g\geq 0$.
	Thus
	\begin{equation}
		0\leq (\sup g)^2\leq \frac{4M^2}{\delta^2}\bigg(
		\frac{4\pi^2}{\Vol(D)}-\frac{\pi}{g-1}+2\pi(2+\Lambda)\sup g\bigg)
	\end{equation}
	Going to the squares, this is a control of the type $X^2\leq A X +B$, with $A,B>0$.
	Explicitely, $A$ and $B$ are the following
	\begin{equation}
		\left\{\begin{array}{ll}
		A&=\frac{8\pi(2+\Lambda)M^2}{\delta^2}\\
		B&=\frac{4M^2}{\delta^1}\bigg(\frac{4\pi^2}{\Vol(D)}-\frac{\pi}{g-1}\bigg)
		\end{array}\right.
	\end{equation}
	so it implies the existence of $C(\delta,\Lambda)=\frac{A+\sqrt{A+4B}}{2}$ such that
	\begin{equation}
		\sup g\leq C(\delta,\Lambda)
	\end{equation}
	Replugging this into the Inequation \ref{eqsupinf}, we get the existence of $\widetilde{C}(\delta,\Lambda)$ satisfying
	\begin{equation}
		|\sup g -\inf g|\leq \widetilde{C}(\delta,\Lambda)
	\end{equation}
	as asserted.
\end{proof}
It happens that the oscillations of $g$ away from $D(z_0,r)$ are
actually the same as those of the log of the section $s$ associated to the line bundle of divisor $\{z_0\}$.
\begin{proposition}\label{fradg}
	%Let $\Sigma$ be a hyperbolic surface with systole $\delta$ and spectral gap $\Lambda$.
	Let $N$ be a degree $1$ line bundle over $\Sigma$ and $s\in H^0(N)$ section, with a
	unique zero $z_0$.
	Denote $D=D(z_0,\delta)$ the disk of radius $\delta$ around $z_0$.
	Denote $a=\frac{2\pi}{\Vol(D)}$.
	Consider $g:\Sigma\rightarrow\R$ the function associated to $z_0$ defined in Proposition \ref{gbdd}.
	%the $C^1$ function satisfying:
	%\begin{equation}
	%	\left\{\begin{array}{cl}
	%	\int_\Sigma g &=0\\
	%	\Delta g &= \frac{1}{2g-2}-a\one_D\end{array}\right.
	%\end{equation}
	Denote  by~$f$ the function $f=-\frac{1}{2}\log|s|^2$.
	Then $f-g$ is bounded and constant on $\Sigma-D(z_0,\delta)$.
\end{proposition}
\begin{proof}
	Let $r$ be the radial geodesic coordinate on $D(z_0,\delta)$.
	With the notations above,
	define the map $h:D\rightarrow\R$:
	\begin{equation}
		h(r)=-2a\ln\cosh(\frac{r}{2})+B+\ln(\tanh\frac{r}{2})\,.
	\end{equation}
	As $h$ is radial, we can explicitly compute the laplacian of $h$:
	\begin{equation}
		\Delta h =\partial^2_r h +\frac{1}{\tanh r}\partial_r h=-a
	\end{equation}
	As we also have
	\begin{equation}
		|z-z_0|=\tanh(\frac{r}{2})
	\end{equation}
	Because $z_0$ is a simple zero of $s$, locally:
	\begin{equation}
		|s(z)|^2\approx K|z-z_0|^2,\quad f(z)=-\ln\tanh{\frac{r}{2}}+\mathcal{O}(1).
	\end{equation}
	And so $f+h$ is bounded.

	Also, we compute the radial derivative at the boundary of the disk of $h$:
	\begin{equation}
		\partial_r h(\delta)=(\frac{1}{2}-a)\tanh\frac{\delta}{2}
		+\frac{1}{2\tanh\frac{\delta}{2}}
	\end{equation}
	But by definition of $a$,
	\begin{equation}
		a=\frac{2\pi}{\Vol(D)}
		=\frac{2\pi}{4\pi\sinh^2(\frac{r}{2})}=\frac{1}{2\sinh^2\frac{r}{2}}
	\end{equation}
	which implies
	\begin{equation}
		\partial_r h(\delta)=0\,.
	\end{equation}
	Hence with the appropriate choice of $B$ we can extend $h$ by zero
	to get a $C^1$ function on $\Sigma$ such that $f+h$ is bounded and satisfies
	\begin{equation}
		\Delta(f+h)=\frac{1}{2g-2}-a\one_D
	\end{equation}
	This implies that $f+h=g+A$ for some constant $:w
	A$.
	As $h$ is constant out of $D$, this means that $f-g$ is constant on $\Sigma-D$,
	as asserted.
\end{proof}
\begin{corollary}\label{oscoutD}
	Let $\Sigma$ be a hyperbolic surface with systole $\delta$ and spectral gap $\Lambda$.
	Let $N$ be a degree $1$ line bundle over $\Sigma$ and $\alpha\in H^0(N)$ section, with a
	unique zero $z_0$.
	Denote $D=D(z_0,\delta)$ the disk of radius $\delta$ around $z_0$.
	Then there is some $C(\delta,\Lambda)$ and $\lambda>0$ such that
	\begin{equation}
		\forall z\in \Sigma-D(z_0,\delta),\quad
		\frac{1}{C(\delta,\Lambda)}\leq\lambda |\alpha|^2\leq C(\delta,\Lambda)
	\end{equation}
\end{corollary}
\begin{proof}
	With the $g$ defined in proposition \ref{gbdd},
	we know
	\begin{equation}
		|\sup g -\inf g|\leq C(\delta,\Lambda)
	\end{equation}
	Denote $f=-\frac{1}{2}\ln|\alpha|^2$.
	By proposition \ref{fradg}, $f-g$ is constant on $\Sigma-D$.
	Thus
	\begin{equation}
		|\underset{\Sigma-D}{\sup} f-\underset{\Sigma-D}{\inf} f|=
		|\underset{\Sigma-D}{\sup} g-\underset{\Sigma-D}{\inf} g|\leq
		C(\delta,\Lambda)
	\end{equation}
	which directly gives, on $\Sigma-D$:
	\begin{equation}
		\frac{\sup |\alpha|^2}{\inf |\alpha|^2}\leq e^{2C(\delta,\Lambda)}
	\end{equation}
	Multiply by $\lambda$ such that $(\sup \lambda|\alpha|^2)(\inf \lambda|\alpha|^2)=1$
	to get the desired result:
	\begin{equation}
		\forall z\in\Sigma-D,\quad
		e^{-C(\delta,\Lambda)}\leq\lambda|\alpha|^2(z)\leq e^{C(\delta,\Lambda)}
	\end{equation}
\end{proof}
It remains to check that we can also have an upper bound on $D$.
We can do this by some version of the Schwarz lemma adapted to sections of positive line bundles:
\begin{proposition}\label{schwarz}
	Let $\Sigma$ be a hyperbolic surface with systole $\delta$ and spectral gap $\Lambda$.
	Let $N$ be a degree $1$ line bundle over $\Sigma$ and $\alpha\in H^0(N)$ section, with a
	unique zero $z_0$.
	Let $D=D(z_0,\delta)$ the disk of radius $\delta$ around $z_0$.
	Then there is $C(\delta)>0$ satisfying:
	\begin{equation}
		\forall z\in D ,\quad
		|\alpha|^2(z)\leq C(\delta)|z|^2\underset{\partial D}{\sup}|\alpha|^2\,.
	\end{equation}
\end{proposition}
\begin{proof}
	Denote $f=-\frac{1}{2}\ln|\alpha|^2$.
	$f$ satisfies
	\begin{equation}
		\Delta f=\frac{1}{2g-2}
	\end{equation}
	Denote $r$ the radial geodesic coordinate on $D$, and consider the map $h$
	\begin{equation}
		h=f+\ln\tanh\frac{r}{2}-\frac{2}{2g-2}\ln\cosh\frac{r}{2}
	\end{equation}
	Then $h$ is bounded on $D$ and satisfies
	\begin{equation}
		\Delta h=0
	\end{equation}
	By the minimum principle,
	\begin{equation}
		\underset{D}{\inf} h =\underset{\partial D}{\inf} h
		=\underset{\partial D}{\inf} f
		+\ln\tanh\frac{\delta}{2}-\frac{1}{g-1}\ln\cosh\frac{\delta}{2}\,.
	\end{equation}
	Exponentiating, we get the upper bound on $|\alpha|$:
	\begin{equation}
		\forall z \in D,\quad
		\frac{|\alpha|^2}{|z|^2}\cosh(\frac{r}{2})^\frac{2}{g-1}\leq
		\underset{\partial D}{\sup}|\alpha|^2
		\frac{(\cosh\frac{\delta}{2})^\frac{2}{g-1}}{(\tanh\frac{\delta}{2})^2}
	\end{equation}
	Because $g>1$ and $\cosh(t)\geq 1$, this implies
	\begin{equation}
		|\alpha|^2(z)\leq |z|^2\underset{\partial D}{\sup}|\alpha|^2
		\frac{(\cosh\frac{\delta}{2})^2}{(\tanh\frac{\delta}{2})^2}
	\end{equation}
\end{proof}
Finally, we combine these results to show our estimate of oscillations of sections
of degree $1$ line bundles:

\begin{proposition}\label{deg1control}
Let $\Sigma$ be a hyperbolic surface. Denote $\delta$ its systole, $\Lambda$ its spectral gap
and consider $r<\frac{\delta}{2}$.
There are  constants $C_i(\delta,\Lambda,r)$ such that, for any holomorphic section $s$
of a degree $1$ line bundle with a unique zero $z_0\in\Sigma$,
there is a $\lambda>0$ such that:
\begin{equation}
|\lambda s|^2\leq C_2\text{ on } B(z_0,r)
\end{equation}
and
\begin{equation}
\frac{1}{C_1}\leq|\lambda s|^2\leq C_1\text{ on }\Sigma-B(z_0,r)
\end{equation}
\end{proposition}
\begin{remark}
Up to renormalizing our section, one can always assume the first condition to be satisfied on
the ball of radius $r$,
and then the results states that the oscillations away from the zero of the section are
controlled, essentially by the spectral gap and the systole of the surface.
\end{remark}
\begin{proof}
Let $\Sigma$ be a hyperbolic surface with systole $\delta$ and spectral gap $\Lambda$.
Let $N$ be a degree $1$ line bundle over $\Sigma$, with $\alpha\in H^0(N)$
a holomorphic section having one unique zero $z_0$.
Consider $D=D(z_0,\delta)$. Thanks to corollary \ref{oscoutD},
up to renormalizing $\alpha$,
there is $C(\delta,\Lambda)>0$ such that
	\begin{equation}
		\forall z\in\Sigma-D,\quad
		\frac{1}{C(\delta,\Lambda)}\leq|\alpha|^2(z)\leq C(\delta,\Lambda)
	\end{equation}
Now because of proposition \ref{schwarz},
	there is $C(\delta)$ such that
	\begin{equation}
		\forall z\in D,\quad
		|\alpha|^2\leq C(\delta)C(\delta,\Lambda)
	\end{equation}
	as asserted.
\end{proof}
Multiplying those sections with covers of a section of a degree $4(g-1)$ line bundle will result
in a balanced section of a bundle of desired degree $4n(g-1)+1$.
\begin{proof}[proof of theorem \ref{equisec}]
	Let $\Sigma$ be a closed hyperbolic surface of genus $g\geq 2$.
	Consider $\Sigma_{n}$ a family of degree $n$ covers of $\Sigma$
	such that the spectral gap of $\Sigma_n$ remains lower bounded by $\Lambda>0$.
	This kind of cover exists, for instance we can use the result of \cite{MNP22},
	which states that for a random cover of $\Sigma$, the probability that
	the relative spectral gap (i.e. the part of the spectral gap which is not in $\Sigma$)
	of the cover is bigger than $\frac{3}{16}-\eps$ goes to $1$ when the degree is large.
	Hence for any $\eta<\min\{\Lambda,\frac{3}{16}\}$, we can choose a family of covers
	whose spectral gap is larger than $\eta$.

	Of course, the systole is nondecreasing by taking a finite cover,
	hence we can consider $\delta$ the systole of $\Sigma$.

	Consider $s$ a 2-form on $\Sigma$, hence a section of $K^2$, which is of degree $4g-4$.
	Denote $s_n$ the lift of $s$ to $\Sigma_n$, which is a section of a degree $4n(g-1)$ line
	bundle.
	Obviously, we have
	\begin{equation}
		\int |s_n|^2=n\int |s|^2,\quad |s_n|_\infty=|s|_\infty
	\end{equation}
	Also, consider $(z_n)$ a sequence of points, $z_n\in\Sigma_n$.
	Define $L_n=\calO(z_n)$ the degree $1$ line bundle associated to the divisor $\{z_n\}$.

	Thanks to Proposition \ref{deg1control}, there are sections $\tau_n$ of $L_n$, and
	constants $C_i(\delta,\Lambda)$ such that
	\begin{equation}
		\left\{\begin{array}{ccc}
			\frac{1}{C_1(\delta,\Lambda)}\leq |\tau_n|^2\leq C_1(\delta,\Lambda)
			&\text{ out of}&D(z_n,\delta)\\
			|\tau_n|^2\leq C_2(\delta,\Lambda)&\text{ on }& D(z_n,\delta)
		\end{array}\right.
	\end{equation}
	Consider now $N_n=K^2L_n$ and $\alpha_n=s_n \tau_n$
	section of $N_n$ line bundle of degree $4n(g-1)+1$.
	By construction, we have the estimates
	\begin{flalign}
		\frac{1}{C_1(\delta,\Lambda)}|s|^2_\infty&\leq|\alpha_n|_\infty^2
		\leq C_2(\delta,\Lambda)|s|^2_\infty\\
		\frac{n-1}{C_1(\delta,\Lambda)}\int|s|^2&\leq\int|\alpha_n|^2
		\leq C_1(n-1)\int|s|^2+C_2\int|s|^2
	\end{flalign}
	Eventually,
	\begin{equation}
		\frac{|\alpha_n|^2_\infty}{\fint|\alpha_n|^2}\leq
		\frac{n}{n-1}\frac{|s|^2_\infty}{\int|s|^2}
		C_1(\delta,\Lambda)C_2(\delta,\Lambda)\Vol(\Sigma)
	\end{equation}
	which is bounded in $n$, so the family $(N_n,\alpha_n)$ is balanced.
\end{proof}
\begin{remark}
	In a broader setting, given a section $s$ of a degree $d$ line bundle,
	we expect the ratio between the $L^2$ average and the $L^\infty$ norm
	to be an quantitative indicator of how much does $s$ look like a lift of a section on
	a smaller surface.
	With this in mind, we expect this ratio to be controlled by
	the lowest distance between the zeroes of $s$.
	This would have the advantage of being an open property,
	and that when it is actually possible, the minimal ratio is gotten by a lift.
	%Considering the properties of the exhibited balanced family of sections,
	%it seems that there should be a stronger result,
	%stating that the ratio between $L^2$ average and $L^\infty$ norm
	%should be bounded by the lowest distance between the zeroes of the section.
	%This would have the advantage of being an open property.
	%However, we couldn't state a precise result in this direction, so we stick with our example.
\end{remark}
\begin{remark}
	The same proof can be adapted to show that for any $a,d\geq 0$,
	we can build a balanced family~$(\Sigma_g,N_g,\alpha_g)$
	with $N_g$ of degree $a(g-1)+d$.
\end{remark}

\section{Estimates on the solutions to the flat curvature equations}
Throughout this section, $\Sigma$ will be a closed hyperbolic surface of genus $g\geq 2$.
Unless otherwise precised, $N$ will be a degree $1$ line bundle over~$\Sigma$ and $\alpha$ will
be a holomorphic $2$-form valued in $N$, i.e. $\alpha\in H^0(K^2N)$.
Fix~$c=\frac{1}{2g-2}$, and consider the following PDE system:
\begin{equation}
\left\{\begin{array}{cl}
\Delta u =& e^{2u}-1+e^{-2u}e^{2v}|\alpha|^2\\
\Delta v =& c-e^{-2u}e^{2v}|\alpha|^2\end{array}\right.
\end{equation}

In this part, we consider each equation separately, considering the first as giving $u$ depending
on $v$, and vice-versa for the second.

\subsection{Dealing with Gauss'equation}

Fixing a map $u$, the first equation is the famous Gauss equation, and is well-behaved.
As the data $v$ and $\alpha$ are fixed, we denote $f=e^{2v}|\alpha|^2$,
and our problem boils down to solve
\begin{equation}
	\Delta u = -1+e^{2u}+e^{-2u}f
\end{equation}
Denote by $C^{0,\alpha}(\Sigma)$ the set of continuous $\alpha$-Hölder functions on $\Sigma$.
Looking at the results of \cite{BS23}, we can apply the same sub-supersolution method to get:
\begin{theorem}\label{Codazzisolve}
	Let $\eta\in(0,1]$. Consider the set
\begin{equation}
	\mathcal{C}:=\{f\in C^{0,\alpha}(\Sigma):0\leq f \leq\frac{\eta}{(1+\eta)^2}\}
\end{equation}
Then for any $f\in\mathcal{C}$, there is a unique $u\in C^{2,\alpha}(\Sigma)$
satisfying:
\begin{flalign}
-\frac{\ln 2}{2}&\leq u \leq 0\\
\Delta u &= e^{2u}-1+e^{-2u}f\label{Gausseq}
\end{flalign}
Moreover, $u$ satisfies the following upper bound:
\begin{equation}
	e^{-2u}f\leq\eta\,.
\end{equation}
\end{theorem}
\begin{proof}
We will proceed with the sub-supersolution method.
A subsolution $v_-$ to the equation~$\Delta u =F(u)$ is a function $v_-$
satisfying $\Delta v_-\geq F(v_-)$.
Analogously, a supersolution is a function $v_+$ satisfying~$\Delta v_+\leq F(v_+)$.

The sub-supersolution principle states that the existence of a subsolution $v_-$ and a
supersolution $v_+$ satisfying $v_-\leq v_+$ ensures the existence of a solution $v$
satisfying $v_-\leq v\leq v_+$.

Freely, $0$ is a supersolution to the equation \ref{Gausseq}.
	As $f\in\mathcal{C}$, the constant map $-\frac{\ln (1+\eta)}{2}$ is a subsolution.
Hence by the sub-supersolution principle, there is a solution $u$
satisfying
\begin{equation}
-\frac{\ln 2}{2}\leq u \leq 0\,.
\end{equation}
To prove uniqueness, assume $u_1$ and $u_2$ are both solutions in $(-\frac{\ln 2}{2},0)$
of
\begin{equation}
\Delta u_i = e^{2u_i}-1+e^{-2u_i}f
\end{equation}
Then the difference $w=u_1-u_2$ satisfies
\begin{equation}
\Delta w = (e^{2u_2}-e^{-2u_2}f)(e^{2w}-1)+fe^{-2u_2}(e^{2w}+e^{-2w}-2)\,.
\end{equation}
We observe that the second term is nonnegative, thus
\begin{equation}
\Delta w \geq (e^{2u_2}-e^{-2u_2}f)(e^{2w}-1)
\end{equation}
but as $f\in\mathcal{C}$ and $u_2\geq\frac{-\ln 2}{2}$, we get
\begin{equation}
e^{-4u_2}f\leq 1
\end{equation}
which implies
\begin{equation}
e^{2u_2}-e^{-2u_2}f|\geq 0
\end{equation}
hence, by the maximum principle, at a maximum of $w$,
\begin{equation}
(e^{2u_2}-e^{-2u_2}f)(e^{2w}-1)\leq 0
\end{equation}
Hence
\begin{equation}
w\leq 0\,.
\end{equation}
Exchanging the roles of $u_1$ and $u_2$ shows
\begin{equation}
u_1=u_2\,.
\end{equation}
So there is a unique solution in $[-\frac{\ln 2}{2},0]$.

	Finally, because $-\frac{\ln(1+\eta)}{2}\leq u\leq 0$,
	we deduce
	\begin{equation}
		e^{-2u}f\leq\frac{\eta}{1+\eta}\leq\eta\,.
	\end{equation}
	As asserted.
\end{proof}
\begin{proposition}\label{Codazzistability}
	The map $\Phi:C^{0,\alpha}(\Sigma,[0,\frac{1}{2}])\rightarrow C^{2,\alpha}(\Sigma)$
	which to $f$ associates the unique $u$ solution of
	\begin{equation}
		\Delta u =e^{2u}-1+e^{-2u}f
	\end{equation}
	satisfying $\frac{-\ln 2}{2}\leq u\leq 0$ is continuous.
\end{proposition}
\begin{proof}
	Thanks to Theorem \ref{Codazzisolve}, $\Phi$ is well-defined. 
	For $f\in C^{0,\alpha}(\Sigma,[0,\frac{1}{2}])$ given, consider the linearized
	equation
	\begin{equation}
		\Delta\dot u =2(e^{2u}-e^{-2u}f)\dot u
	\end{equation}
	The assumptions on $u$ imply 
	\begin{equation}
		e^{2u}-e^{-2u}f\geq 0
	\end{equation}
	By a maximum principle, the only solution of the linearized equation is $0$,
	hence our solutions are stable and the defined map is continuous.
\end{proof}
\subsection{Dealing with Ricci's equation}
Considering the second equation, we will fix $v$ and express a solution $v$ as the maximizer of
the following functional:
\begin{notation}
For $u\in C^{1,\alpha}(\Sigma)$, consider $W^{1,2}_0(\Sigma)$ the Sobolev space of zero average
functions that are square integrable with square integrable gradient.
Define the functional $J_u:W^{1,2}_0(\Sigma)\rightarrow\R$ as
\begin{equation}
J_u(w):=\log\big(\fint e^{-2u}|\alpha|^2e^{2w}\big)-\frac{1}{c\Vol\Sigma}\int|\nabla w|^2
\end{equation}
with $c=\frac{1}{2g-2}$.
The notation $\fint$ is used for the average of the considered function.

We will also use the notation $|w|_2^2$ for the squared $L^2$-norm:
\begin{equation}
	|w|_2^2=\int_\Sigma |w|^2
\end{equation}
\end{notation}
When no confusion is possible, we will write $J$ for $J_u$.
We will prove that $J$ attains its maximum at a function $w$ which will be, up to an addittive
constant, the solution $v$ that we are looking for.
This proof is heavily inspired from the method used to solve the prescribed curvature equation
on the sphere when the prescribed data is even. Chang--Yang \cite{CY03} use this method in this
precise case. Keep in mind that our equation is not a prescribed curvature equation on the
surface, but rather a prescribed curvature equation on a disk bundle on it.

\begin{lemma}\label{JweakC0}
$J$ is well-defined and weakly continuous on $W^{1,2}_0(\Sigma)$.
\end{lemma}
\begin{proof}
First, we need to check that $J$ is well-defined.
This is true as $\Sigma$ is compact, and the Sobolev constants in the injections
$W^{1,2}\rightarrow L^p$ are exponentially integrable (see Ilias \cite{Ili83})
\begin{equation}\label{expint}
\exists C>0: \forall w\in W^{1,2}_0,\text{ with }\int |\nabla w|^2\leq 1:\int e^{2w}\leq C
\end{equation}
As $e^{-2u}|\alpha|^2$ is bounded on $\Sigma$, we deduce that $J(w)$ is finite.
Hence $J$ is well-defined.
It is standard that the $L^2$~norm of the gradient is weakly continuous on $W^{1,2}(\Sigma)$.
To be sure that~$J$ is weakly continous, it remains to prove that the map
\begin{equation}
H:w\mapsto\log\big(\fint e^{-2u}|\alpha|^2e^{2w}\big)
\end{equation}
is.
A simple computation shows that $H$ is actually convex, hence it is weakly continuous if and
only if it is strongly continuous.
It is strongly continuous because of the aforementioned property, equation (\ref{expint}).
\end{proof}
In order to study explicitly the behavior of $J$, we will need a more precise control on
the Sobolev constants of the functions involved, namely we will use a precise Moser-Trudinger
inequality obtained in \cite{Bro23}.
\begin{theorem}[Moser--Trudinger, B.]\label{MT}
Denote respectively by $\delta$ and $\Lambda$ the systole and spectral gap of $\Sigma$.
There is a constant $C(\delta,\Lambda)$ depending only on $\delta, \Lambda$ and continuously in
those parameters
such that:
\begin{equation}
\forall u\in W^{1,2}_0(\Sigma),\,\int|\nabla u|^2\leq 1\Rightarrow\int (e^{4\pi u^2}-1)\leq C(\delta,\Lambda)
\end{equation}
\end{theorem}
\begin{corollary}\label{Jbdd}
%Assume $N$ is of degree $1$, and $c=\frac{1}{2g-2}$.
The functional $J$ is upper bounded on $W^{1,2}_0(\Sigma)$, and
we have the following estimate:
\begin{equation}
	\log\big(\fint e^{-2u}|\alpha|^2\big)\leq\sup J\leq \log\max e^{-2u}|\alpha|^2
	+\log\big(1+\frac{C(\delta,\Lambda)}{\Vol(\Sigma)}\big)\,.
\end{equation}
	Also, $J$ is proper on $W_0^{1,2}(\Sigma)$, as the following holds
	\begin{equation}
		\forall w \in W_0^{1,2}(\Sigma),\, J(w)+\frac{1}{4\pi}|\nabla w|_2^2\leq
		\log\max e^{-2u}|\alpha|^2
	+\log\big(1+\frac{C(\delta,\Lambda)}{\Vol(\Sigma)}\big)\,.
	\end{equation}
\end{corollary}
\begin{proof}
Let $w\in W^{1,2}_0(\Sigma)$.
Consider that
\begin{equation}
2w\leq 4\pi\frac{w^2}{|\nabla w|_2^2}+\frac{|\nabla w|_2^2}{4\pi}
\end{equation}
Exponentiating this inequality, we get
\begin{equation}
	e^{2w}\leq e^{\frac{|\nabla w|_2^2}{4\pi}} e^{4\pi\big(\frac{w}{|\nabla w|_2}\big)^2}
\end{equation}
	Using theorem \ref{MT} to bound the integral of the right-hand side,
\begin{equation}
\int e^{2w}\leq (C(\delta,\Lambda)+\Vol(\Sigma))e^\frac{|\nabla w|_2^2}{4\pi}
\end{equation}
We deduce
\begin{equation}
J(w)\leq\log\max e^{-2u}|\alpha|^2
	+\big(\frac{1}{4\pi}-\frac{1}{c\Vol\Sigma}\big)|\nabla w|_2^2
	+\log\big(1+\frac{C(\delta,\Lambda)}{\Vol\Sigma}\big)
\end{equation}
	As $\Vol(\Sigma)= 4\pi(g-1)$,and the bundle is of degree $1$, necessarily
	\begin{equation}
		\frac{1}{4\pi}-\frac{1}{c\Vol(\Sigma)}=-\frac{1}{4\pi}
	\end{equation}
	Inserting it on the left-hand side, we get
	\begin{equation}
	J(w)+\frac{1}{4\pi}|\nabla w|_2^2\leq
		\log\max e^{-2u}|\alpha|^2
	+\log\big(1+\frac{C(\delta,\Lambda)}{\Vol(\Sigma)}\big)\,.
	\end{equation}
	So $J$ is indeed proper and upper-bounded with the desired explicit upper bound.

The lower estimate for $\sup J$ is gotten by expliciting $J(0)\leq \sup J$.
\end{proof}
\begin{proposition}\label{Jmax}
	%Assume $N$ is of degree $1$, and $c=\frac{1}{2g-2}$.
	The functional $J$ attains its maximum at $w\in W^{1,2}_0(\Sigma)$ satisfying
	\begin{equation}
		\frac{1}{4\pi}|\nabla w|_2^2\leq
		\log\frac{\max e^{-2u}|\alpha|^2}{\fint e^{-2u}|\alpha|^2}
		+\log\big(1+\frac{C(\delta,\Lambda)}{\Vol\Sigma}\big)\,.
	\end{equation}
\end{proposition}
\begin{proof}
	Thanks to Corollary \ref{Jbdd}, we know that $J$ is upper bounded.
	Consider a maximizing sequence $(w_n)$ of $J$ is $W^{1,2}_0(\Sigma)$.
	Writing down our estimate
	\begin{equation}
		J(w_n)+\frac{1}{4\pi}|\nabla w_n|_2^2\leq\log\max e^{-2u}|\alpha|^2
		+\log\big(1+\frac{C(\delta,\Lambda)}{\Vol\Sigma}\big)
	\end{equation}
	As $(J(w_n))$ is upper bounded by $\sup J$, so is $(|\nabla w_n|_2^2$,
	We deduce that  $(|\nabla w_n|_2^2)$ is bounded.
	The explicit bound, assuming $J(w_n)\geq J(0)$, is
	\begin{equation}
		\frac{1}{4\pi}|\nabla w_n|_2^2\leq
		\log\frac{\max e^{-2u}|\alpha|^2}{\fint e^{-2u}|\alpha|^2}
		+\log\big(1+\frac{C(\delta,\Lambda)}{\Vol\Sigma}\big)\,.
	\end{equation}
	As $(w_n)$ is bounded in $W_0^{1,2}(\Sigma)$,
	extracting a subsequence we can assume it converges weakly
	towards some~$w\in W^{1,2}_0(\Sigma)$.
	By lemma \ref{JweakC0}, $J$ is weakly continuous, so $w$ is a maximizer of~$J$.
	As $|\nabla w|_2^2\leq \liminf |\nabla w_n|_2^2$, we get the desired bound on $w$.
\end{proof}
\begin{theorem}\label{Ricciexistence}
	Let $N$ be of degree $1$, $\alpha\in H^0(K^2 N)$, $c=\frac{1}{2g-2}$
	and $u\in C^{0,\alpha}(\Sigma, [\frac{-\ln 2}{2},0])$.
	Then
	there is a solution $v$
	to
	\begin{equation}
		\Delta v = c-e^{2v}e^{-2u}|\alpha|^2
	\end{equation}
	satisfying
	\begin{equation}
	\frac{1}{4\pi}|\nabla v|_2^2\leq
		\log\frac{\max e^{-2u}|\alpha|^2}{\fint e^{-2u}|\alpha|^2}
		+\log\big(1+\frac{C(\delta,\Lambda)}{\Vol\Sigma}\big)\,.
	\end{equation}
\end{theorem}
\begin{proof}
	Thanks to proposition \ref{Jmax},
	there exist $w\in W^{1,2}_0(\Sigma)$ a maximizer of the functional $J$.
	For any $s\in W^{1,2}_0$, $D_wJ\cdot s=0$ and so
	\begin{equation}
		\frac{2}{\Vol(\Sigma)}
		\frac{\int e^{-2u}|\alpha|^2e^{2w}s}{\fint e^{-2u}|\alpha|^2e^{2w}}
		-\frac{2}{c\Vol(\Sigma)}\int \nabla w\cdot\nabla s=0
	\end{equation}
	Weakly, this means
	\begin{equation}
		\frac{1}{c}\Delta w+\frac{e^{-2u}|\alpha|^2e^{2w}}{\fint e^{-2u}|\alpha|^2e^{2w}}\in (L^2_0)^\perp 
	\end{equation}
	As the orthogonal of the zero average functions are the constant maps,
	there is a constant~$A$ such that
	\begin{equation}
		\Delta w =A-c\frac{e^{-2u}|\alpha|^2e^{2w}}{\fint e^{-2u}|\alpha|^2e^{2w}}
	\end{equation}
	By the elliptic regularity principle, as $u$ is in $C^{0,\alpha}(\Sigma)$,
	$w$ belongs to $C^{2,\alpha}(\Sigma)$ and satisfies the above equation strongly.

	Writing $\int \Delta w =0$, we deduce
	\begin{equation}
		A=c\,.
	\end{equation}
	Consider the element of $W^{1,2}(\Sigma)$:
	\begin{equation}
		v=w+\frac{1}{2}\ln\big(\frac{c}{\fint e^{-2u}|\alpha|^2e^{2w}}\big)
	\end{equation}
	As $\Delta v=\Delta w$, we get
	\begin{equation}
		\Delta v =c-e^{-2u}|\alpha|^2e^{2v}
	\end{equation}
	As $|\nabla w|_2^2=|\nabla v|_2^2$, the estimate of Proposition \ref{Jmax}
	ensures the desired bound on $|\nabla v|_2^2$.
\end{proof}
\begin{remark}
	In contrast with the Gauss equation, we cannot in this case establish the
	uniqueness of the solution $v$. We will address this problem
	by showing stability of the solutions under some conditions after stating some 
	precise estimates on the solution $w$ and~$v$.
\end{remark}
\subsection{Estimating the regularity of the solution}
In order to apply later a fixed point theorem,
we need to understand the behavior of the solutions $v$.
Particularly, we need an $L^\infty$ bound on $e^{2v}|\alpha|^2$.
Because of the sign difference with Gauss'equation,
applying a maximum principle won't yield bounds on $v$.
We get our upper bound via a $W^{2,2}$-estimate on the maximizer $w$ of $J_u$.

\begin{proposition}\label{J22}
%	Let $N$ be of degree $1$,
%	and $c=\frac{1}{2g-2}$.
	Assume $u\in C^{0,\alpha}(\Sigma)$. Then a maximizer $w$
	of the functional $J_u$ satisfies
	\begin{equation}
		|\nabla^2 w|_2^2\leq
		c^2\bigg((\Vol\Sigma+C(\delta,\Lambda))(\frac{\max|\alpha|^2}{\fint |\alpha|^2})^2-\Vol\Sigma
		\bigg)+4\pi\bigg(
		\log\frac{\max e^{-2u}|\alpha|^2}{\fint e^{-2u}|\alpha|^2}
		+\log\big(1+\frac{C(\delta,\Lambda)}{\Vol\Sigma}\big)\bigg)
	\end{equation}
\end{proposition}
\begin{proof}
	Let $w$ be a maximizer of $J$.
	First, we use the Bochner identity to estimate the Hessian of $w$
	\begin{equation}
		|\nabla^2w|_2^2=\int(\Delta w)^2+\int |\nabla w|^2\,.
	\end{equation}
	Since $\Delta w$ satisfies
	\begin{equation}
		\Delta w =c-\frac{c|\alpha|^2 e^{-2u}e^{2w}}{\fint e^{-2u}|\alpha|^2 e^{2w}}\,.
	\end{equation}
	We deduce that
	\begin{equation}
		\int (\Delta w)^2= -c^2\Vol\Sigma+c^2(\Vol\Sigma)^2\frac{\int e^{-4u}|\alpha|^4e^{4w}}{(\int e^{-2u}|\alpha|^2e^{2w})^2}
	\end{equation}
	To control the right hand side, we  estimate $4w$
	\begin{equation}
		4w\leq\frac{4\pi w^2}{|\nabla w|_2^2}+\frac{|\nabla w|_2^2}{\pi}
	\end{equation}
	and apply the Moser--Trudinger inequality of Theorem \ref{MT} to get
	\begin{equation}
		\int e^{4w}\leq (\Vol\Sigma+C(\delta,\Lambda))e^{\frac{|\nabla w|_2^2}{\pi}}
	\end{equation}
	This ensures
	\begin{equation}
		\int |\alpha|^4e^{-4u}e^{4w}\leq(\Vol\Sigma+C(\delta,\Lambda))\max e^{-4u}|\alpha|^4\exp\big(\frac{|\nabla w|_2^2}{\pi}\big)\,.
	\end{equation}
	It remains to estimate the denominator of the right hand side.
	As $c\Vol\Sigma=2\pi$ and by maximality $J(w)\geq J(0)$, we have
	\begin{equation}
		\log\fint|\alpha|^2 e^{-2u}\leq \log\fint e^{-2u}|\alpha|^2e^{2w}-\frac{1}{2\pi}|\nabla w|_2^2\,.
	\end{equation}
	Exponentiating this,
	\begin{equation}
		\int e^{-2u}|\alpha|^2\exp\big(\frac{1}{2\pi}|\nabla w|_2^2\big)\leq\int e^{-2u}|\alpha|^2e^{2w}
	\end{equation}
	Combining this with the upper bound on the numerator, we get
	\begin{equation}
		\frac{\int e^{-4u}|\alpha|^4e^{4w}}{\big(\int e^{-2u}|\alpha|^2e^{2w}\big)^2}
		\leq
		\frac{(\Vol\Sigma+C(\delta,\Lambda))\max e^{-4u}|\alpha|^4\exp(\frac{|\nabla w|_2^2}{\pi})}
		{\big(\int e^{-2u}|\alpha|^2\exp(\frac{1}{2\pi}|\nabla w|_2^2)\big)^2}
	\end{equation}
	It happens that the left-hand side simplifies to give
	\begin{equation}
		\int (\Delta w)^2\leq -c^2\Vol\Sigma+c^2(\Vol\Sigma+C(\delta,\Lambda))
		\frac{\big(\max e^{-2u}|\alpha|^2\big)^2}{(\fint e^{-2u}|\alpha|^2)^2}
	\end{equation}
	Combining with the bound on $|\nabla w|_2^2$ from Proposition \ref{Jmax},
	we get the desired result.
\end{proof}
We sum up with the local upper bound on the solution,
which will be key to building solutions
of the PDE system:
\begin{proposition}\label{pointwisecontrol}
%Let $\Sigma$ be a closed hyperbolic surface of genus $g$, $N$ a line bundle over $\Sigma$
%endowed with its constant curvature metric,
%and $\alpha$ a section of $N$. Denote $c=\frac{1}{2g-2}$.
	Let $u:\Sigma\rightarrow\R$ be a map in $C^{0,\alpha}(\Sigma)$ such that
\begin{equation}
	1\leq e^{-2u}\leq 2\,.
\end{equation}
Then there is a constant $C_1=C_1(\delta,\Lambda,\frac{|\alpha|_\infty^2}{\fint|\alpha|^2})$, 
and $v\in C^{2,\alpha}(\Sigma)$ a solution of
\begin{equation}
	\Delta v = c-e^{-2u}e^{2v}|\alpha|^2
\end{equation}
satisfying
\begin{equation}
	\sup e^{2v}|\alpha|^2\leq \frac{C_1(\delta,\Lambda,\frac{|\alpha|_\infty^2}{\fint|\alpha|^2})}{2g-2}
\end{equation}
	The constant $C_1$ is continuous in the parameters $\delta,\Lambda$ and
	$\frac{|\alpha|_\infty^2}{\fint|\alpha|^2}$.
\end{proposition}
\begin{proof}
	Once again, consider $w\in W_0^{1,2}$ a maximizer of $J$
	And build the solution $v$ as a translated of $w$
	so that
	\begin{equation}
		w=v-\overline v\,.
	\end{equation}
	Thanks to Proposition \ref{J22},
	there is a bound $C_1(\delta,\Lambda,\frac{|\alpha|_\infty^2}{\fint|\alpha|^2})$
	such that
	\begin{equation}
		|w|_{2,2}\leq C_1(\delta,\Lambda,\frac{|\alpha|_\infty^2}{\fint|\alpha|^2})
	\end{equation}
	As the continuous embedding $W^{2,2}\hookrightarrow L^\infty$ depends only on the systole
	$\delta$, this implies that we have $C_2>0$ with the same dependencies such that
	\begin{equation}
		|w|_\infty\leq C_2(\delta,\Lambda,\frac{|\alpha|_\infty^2}{\fint|\alpha|^2})
\,.
	\end{equation}
	We compute the average of $e^{2v}|\alpha|^2$
	\begin{equation}
		c=\fint e^{-2u}e^{2v}|\alpha|^2=e^{2\overline v}\fint e^{-2u}e^{2w}|\alpha|^2
	\end{equation}
	from which we get the estimate:
	\begin{equation}
	e^{2\overline v}e^{-2C_2}\fint|\alpha|^2\leq\frac{1}{2g-2}\leq
		2e^{2\overline v}e^{2C_2}\fint|\alpha|^2
	\end{equation}
	Combining this, we get
	\begin{equation}
		e^{2v}|\alpha|^2 = e^{2\overline v}e^{2w}|\alpha|^2\leq
		\frac{e^{4C_2}}{2g-2}\frac{|\alpha|^2_\infty}{\fint|\alpha|^2}
	\end{equation}
	as asserted.
	As the constant $C_1$ from Proposition \ref{J22} is continuous in the parameters,
	so is the exhibited constant $C_2$.
\end{proof}
It remains to show the stability of the solutions under suitable perturbations.
\begin{lemma}\label{Riccistability}
	%Denote $\Lambda$ the spectral gap of $\Sigma$.
	Consider $f\geq 0$ a function in $C^{0,\alpha}(\Sigma)$,
	such that there is $v$ solution of
	\begin{equation}
		\Delta v = c-e^{2v}f
	\end{equation}
	satisfying $2e^{2v}f<\Lambda$.
	Consider the operator $H=\Delta+2e^{2v}f:W^{2,2}(\Sigma)\rightarrow L^2(\Sigma)$.
	Then $H$ is invertible, and we have the upper bound on its inverse:
	\begin{equation}
		||H^{-1}||\leq\max\{(\Lambda-2|e^{2v}f|_\infty)^{-1},c^{-1}\}
	\end{equation}
	%Then there are open sets $\mathcal{U},\mathcal{V}$
	%where $\mathcal{U}$ is a neighborhood of $f$ in $C^{0,\alpha}(\Sigma,\R_+)$
	%and $\mathcal{V}$ neighborhood of $v$ in $C^{2,\alpha}(\Sigma)$
	%and a continuous map $\Psi:\mathcal{U}\rightarrow\mathcal{V}$
	%such that for any $\widetilde f\in\mathcal{U}$, $\Psi(\widetilde f)=\widetilde v$
	%satisfies:
	%\begin{equation}
	%	\left\{\begin{array}{lc}
%			\Delta\widetilde v&= c-e^{2\widetilde v}\widetilde f\\
%			2e^{2\widetilde v}\widetilde f&<\Lambda
%		\end{array}\right.
%	\end{equation}
	%Don't lift now to a continuous branch, keep that for part 5 !
\end{lemma}
\begin{proof}
	$H$ is a self-adjoint operator between Hilbert spaces.
	Moreover, for $\lambda$ large enough,
	we have the estimate:
	\begin{flalign}
		|(\lambda I-H)w|_2^2 &|(\lambda-e^{2v}f)w-\Delta w|_2^2\\
		&\leq (\lambda-\frac{\Lambda}{2})^2\int w^2
		+\int (\Delta w)^2-2\int(\lambda-e^{2v}f)w\Delta w\\
		&\leq C|w|_{2,2}^2
	\end{flalign}
	Hence for $\lambda$ large enough, $(\lambda I-H)$ has compact inverse,
	so it is a Fredholm operator, and $H$ has discrete spectrum.
	
	To conclude the lemma, it remains to check that $H$ has no eigenvalue
	in the range $(-\Lambda+|2e^{2v}f|_\infty,2c)$.

	Let $t$ be in this range. If $\dot v$ is an eigenvector for the value 
	$t$,
	we have
	\begin{equation}
		\Delta\dot v+ (2e^{2v}f-t)\dot v=0
	\end{equation}
	Decompose $\dot v =m(\dot v)+ \check{v}$ with $m(\dot v)$ being the average of $\dot v$.
	Then it comes
	\begin{flalign}
		m(\dot v)\int 2e^{2v}f-t&=-\int (2e^{2v}f-t)\check{v}\\
		\int |\nabla \check{v}|^2-(2e^{2v}f-t)\check{v}^2&=
		m(\dot v)\int (2e^{2v}f-t)\check{v}
	\end{flalign}
	Combining these two equations with $\int e^{2v}f=c$, we get
	\begin{equation}
		\int |\nabla \check{v}|^2-(2e^{2v}f-t)\check{v}^2=-2m(v)^2(2c-t)\Vol\Sigma
	\end{equation}
	As $t<2c$, the right-hand side is nonpositive,
	while because $-\Lambda+2e^{2v}f<t$, the left-hand side is nonnegative.
	So both must be zero, and $\dot v=0$, as desired.
\end{proof}
A consequence of this lemma is that we locally able to lift solutions
to $\Delta v= c-e^{2v}f$. But we will need to be able to lift to a continuous branch for all $f$
obtained from solving the Gauss-Codazzi equation. To this aim, we give a more
precise description of the open set on which we are able to lift:
\begin{theorem}\label{globallift}
	Let $\Sigma$ be a closed hyperbolic surface. We will denote by $(C_i)$ constants
	depending on the systole, the spectral gap, and the ratio
	$\frac{|\alpha|^2_\infty}{\fint|\alpha|^2}$.
	Consider $v_0$ a solution of
	\begin{equation}
		\Delta v_0 = c-e^{2v_0}|\alpha|^2
	\end{equation}
	satisfying
	\begin{equation}
		e^{2v}|\alpha|^2\leq \frac{C_1}{2g-2}<\Lambda\,.
	\end{equation}
	Fix $\eps>0$.
	Consider the functional $\calH:C^{0,\alpha}(\Sigma)\times W_{2,2}(\Sigma)\rightarrow L^2(\Sigma)$:
	\begin{equation}
		\calH(f,v)=\Delta v -c+e^{2v}f|\alpha|^2
	\end{equation}
	Then there is $C_2>0$,
	such that for $|f-1|_\infty\leq\frac{C_2}{\sqrt{2g-2}}$,
	there is a unique $v\in W^{2,2}$ satisfying:
	\begin{equation}
		\left\{\begin{array}{ll}
			|v-v_0|_{2,2}&\leq\eps\\
			\Delta v&=c-e^{2v}f|\alpha|^2
		\end{array}\right.
	\end{equation}
		In particular, this $v$ also satisfies:
	\begin{equation}
		e^{2v}|\alpha|^2\leq \frac{C_3}{2g-2}
	\end{equation}
\end{theorem}
\begin{proof}
	This explicit lift comes from an estimation of the bounds appearing in the implicit function
	theorem. We will use the version in \cite{Hol70}.
	With the notations of the theorem,
	we need to compute $k_1=|(D_v\calH(f_0,v_0))^{-1}|$.
	The lemma \ref{Riccistability} shows that,
	provided~$\frac{1}{2(2g-2)}\leq\Lambda$, 
	\begin{equation}
		k_1\leq\frac{1}{2}(2g-2)\,.
	\end{equation}
	Now we need to find $\delta>0$ and $g_1$ such that, for $s\in[0,\delta]$, $v\in[0,\eps]$

	\begin{equation}
		|f-1|\leq\delta,\,|v-v_0|_{2,2}\leq\eps\Rightarrow
		||D_v\calH(f,v)-D_v\calH(1,v_0)||\leq g_1(|f-1|_\infty,|v-v_0|_{2,2})
	\end{equation}
	But, as
	\begin{equation}
		|D_v\calH(f,v)\dot v-D_v\calH(1,v_0)\dot v|_2^2=
		\int (2(e^{2v}f-e^{2v_0})|\alpha|^2\dot v)^2
		\leq
		\frac{4C_1^2}{(2g-2)^2}e^{4|v-v_0|_\infty}|1-f|_\infty|\dot v|_2^2
	\end{equation}
	from which we deduce
	\begin{equation}
		||D_v\calH(f,v)-D_v\calH(1,v_0)||\leq \frac{C_4}{2g-2}e^{2C_5|v-v_0|_{2,2}}|1-f|_\infty
	\end{equation}
	Finally, we must find $g_2$ such that
	\begin{equation}
		|f-1|_\infty\leq\delta\Rightarrow ||\calH(f,v_0)||\leq g_2(|f-1|_\infty)
	\end{equation}
	Explicitly,
	\begin{equation}
	|\calH(f,v_0)|_2^2=\int (e^{2v_0}|\alpha|^2(f-1)^2)^2\leq \frac{C_6}{2g-2}
	\end{equation}
	Of which we get
	\begin{equation}
		g_2(s)=\frac{C_7s}{\sqrt{2g-2}}
	\end{equation}
	Now we can lift a continuous branch satisfying $\calH(f,\Psi(f))=0$
	provided $\delta$ satisfies
	\begin{equation}
		k_1g_1(\delta,\eps)\leq\eta<1\text{ and }k_1g_2(\delta)\leq\eps(1-\eta)
	\end{equation}
	This is verified for $\delta$ of the kind:
	\begin{equation}
		\delta =\frac{C_8}{\sqrt{2g-2}}
	\end{equation}
	and the solutions $v=\Psi(f)$
	satisfy
	\begin{equation}
		|v-v_0|_{2,2}\leq\eps
	\end{equation}
	Of which we deduce
	\begin{equation}
		e^{2v}|\alpha|^2\leq e^{2C_9|v-v_0|_{2,2}}e^{2v_0}|\alpha|^2
		\leq \frac{C_{10}}{2g-2}
	\end{equation}
	as desired.
\end{proof}

\section{Proof of the main theorem}
In this section we will prove the geometric meaning of all this study:
\begin{theorem}\label{analysisthm}
	Let $\Sigma$ be a closed hyperbolic surface of genus $g$. Fix $\eta\in (0,1)$ and consider a family
$(\Sigma_n, N_n,\alpha_n)_{n\geq 2}$ satisfying:
\begin{itemize}
\item[(a)]
$\Sigma_n$ is a degree $n$ cover of $\Sigma$, equipped with the lift of
the hyperbolic metric on $\Sigma$.
\item[(b)]
$N_n$ is a degree $4n(g-1)+1$ line bundle over $\Sigma_n$, equipped with a metric of constant
sectional curvature.
\item[(c)]
$\alpha_n\in H^0(N_n)$ is a section of $N_n$.
\item[(d)]
There is $\Lambda_0>0$ such that for all $n$ the spectral gap of $\Sigma_n$ is bigger than $\Lambda_0$.
\item[(e)]
The family of sections is balanced, i.e.
\begin{equation}
\underset{n\geq 2}{\sup}\frac{|\alpha_n|^2_\infty}{\fint|\alpha_n|^2}<\infty\,.
\end{equation}
\end{itemize}
Then, for $n$ large enough, there are $u,v:\Sigma_n\rightarrow \R$ such that
\begin{equation}
\underset{z\in\Sigma_n}{\sup}\,e^{-4u}e^{2v}|\alpha|^2(z)<\eta
\end{equation}
\begin{equation}
\left\{\begin{array}{cl}
\Delta u&=e^{2u}-1+e^{-2u}e^{2v}|\alpha_n|^2\\
\Delta v &=\frac{1}{n(2g-2)}-e^{-2u}e^{2v}|\alpha_n|^2\end{array}\right.
\end{equation}
\end{theorem}
As a corollary, we get the desired existence theorem for almost-fuchsian representations with
nontrivial normal bundle.
\begin{corollary}
There is a genus $g_0$ such that for any $g\geq g_0$,
there exists a representation $\rho:\pi_1\Sigma_g\rightarrow SO(4,1)$
which is almost-fuchsian, a complex variation of Hodge structure and the hyperbolic $4$-manifold
$\rho\backslash\H^4$ is a degree $1$ disk bundle over $\Sigma_g$.
\end{corollary}
\begin{proof}
	Assuming Theorem \ref{analysisthm}, let $\Sigma$ be a closed hyperbolic surface of genus $2$.
Thanks to the result of Magee--Naud--Petri \cite{MNP22},
we know the existence of a family of covers~$(\Sigma_n)$ satisfying condition $(d)$ of
Theorem \ref{analysisthm}.
Thanks to Theorem \ref{equisec}, there is a balanced sequence~$(\alpha_n,N_n)$
satisfying conditions (b), (c) and (e).
Thus for $n$ large enough we have~$(u,v)$ satisfying
\begin{equation}
\underset{z\in\Sigma_n}{\sup}\,e^{-4u}e^{2v}|\alpha|^2(z)<1
\end{equation}
and
\begin{equation}
\left\{\begin{array}{cl}
\Delta u&=e^{2u}-1+e^{-2u}e^{2v}|\alpha_n|^2\\
\Delta v &=\frac{1}{n(2g-2)}-e^{-2u}e^{2v}|\alpha_n|^2\end{array}\right.
\end{equation}
Thanks to Theorem \ref{fundathm}, this corresponds to an almost-fuchsian representation and
complex  variation of Hodge structure, with holomorphic data $\alpha$, and the hyperbolic
$4$-manifold is a degree $1$ disk bundle over $\Sigma_n$.
\end{proof}
The main argument of Theorem \ref{analysisthm} will be a fixed point theorem, namely
Schauder's fixed point theorem:
a continuous map which preserves a compact convex set has a fixed point.
\begin{proof}[ Proof of Theorem \ref{analysisthm} ]
	Denote $\mathcal{K}^n_u$ the following convex set:
	\begin{equation}
		\mathcal{K}^n_u:=\{u\in C^{0,\alpha}(\Sigma_n):-\frac{\ln (1+\eta)}{2}\leq u\leq 0\}
	\end{equation}
	Thanks to Theorem \ref{Ricciexistence}, for $u=0$,
	there is a function $v_0$ solution of
	\begin{equation}
		\Delta v_0 = \frac{1}{|\chi(\Sigma_n)|}-e^{2v_0}|\alpha_n|^2
	\end{equation}
	and by Proposition \ref{pointwisecontrol}, there is $C_0>0$ depending continuously on
	$\delta,\Lambda$ and the ratio $\frac{|\alpha_n|_\infty^2}{\fint|\alpha_n|^2}$
	such that
	\begin{equation}
		|e^{2v_0}|\alpha_n|^2\leq \frac{C_0}{2n(g-1)}
	\end{equation}
	The assumptions $(a)$, $(d)$ and $(e)$ ensure that the constant $C$ can be chosen
	to be independent of $n$, as it must be for our proof.
	Recall $\Lambda_0$ defined in $(d)$ as a common lower bound to the spectral gaps of $(\Sigma_n)$.
	Consider now $n$ large enough so that:
	\begin{equation}
		\frac{C_0}{2n(g-1)}<\min\{\Lambda_0,\frac{1}{2}\}
	\end{equation}
	%This inequality ensures, by lemma \ref{Riccistability},
	%that the solutions $v$ are stable, for any $u\in\mathcal{K}^n_u$
	%and be locally continuously lifted.
	%But as $\mathcal{K}^n_u$ is convex, it is contractible, and we can extend globally our continuous lifts to a continuous map $\Psi:\mathcal{K}^n_u\rightarrow C^{2,\alpha}(\Sigma_n)$.

	Thanks to proposition \ref{globallift}, there are constants $C_1,C_2,C_3$ independent of $n$,
	and a continuous map
	 \begin{equation}
		 \Psi:\{u\in C^{0,\alpha}||e^{2u}-1|_\infty\leq \frac{C_1}{\sqrt{2n(g-1)}}\}
		 \rightarrow\{v||v-v_0|_\infty\leq \frac{C_2}{2n(g-1)}\}
	 \end{equation}
		 such that $v=\Psi(u)$ is a solution of
	\begin{equation}
		\Delta v =\frac{1}{|\chi(\Sigma_n)|}-e^{2v}e^{-2u}|\alpha_n|^2
	\end{equation}
	And $v=\Psi(u)$ satisfies
	\begin{equation}
		|e^{2v}|\alpha_n|^2\leq\frac{C_3}{2n(g-1)}
	\end{equation}
	Provided $\frac{C_3}{2n(g-1}\leq\frac{\eta}{(1+\eta)^2}$, theorem \ref{Codazzisolve}
	gives us a map $\Phi$ and an $n$-independent constant $C_4$
	\begin{equation}
		\Phi:\{v|e^{2v}|\alpha_n|^2\leq\frac{C_3}{2n(g-1)}\}\rightarrow
		\{u||e^{2u}-1|_\infty\leq\frac{C_4}{2n(g-1)}\}
	\end{equation}
	such that $u=\Phi(v)$ is a solution of
	\begin{equation}
		\Delta u = e^{2u}-1+e^{-2u}e^{2v}|\alpha_n|^2
	\end{equation}
	and $u$ satisfies
	\begin{equation}
		\sup e^{-4u}e^{2v}|\alpha_n|^2<\eta
	\end{equation}
	Finally for $n$ sufficiently large,
	\begin{equation}
		\frac{C_4}{2n(g-1)}\leq\frac{C_1}{\sqrt{2n(g-1)}}
	\end{equation}
	and $\Phi\circ\Psi$ is a well-defined self-map of $\mathcal K=\{u| |e^{2u}-1|\leq\frac{C_4}{2n(g-1)}\}\subset \mathcal K_n^u$.
	We also point out that ,by ellipticity, every element in the image of $\Phi$
	is in $C^{2,\alpha}$ and satisfies
	\begin{equation}\left\{\begin{array}{lc}
		-\frac{\ln 2}{2}&\leq u\leq 0\\
	-1&\leq \Delta u\leq 1\end{array}\right.
	\end{equation}
	And so the image of $\Phi$ is in the following set :
	\begin{equation}
		\mathcal{Q}=\{u\in C^{2,\alpha}(\Sigma_n): -1\leq\Delta u\leq 1\text{ and }
		-\frac{\ln 2}{2}\leq u\leq 0\}
	\end{equation}
	whose closure in $C^{0,\alpha}(\Sigma_n)$ is a convex compact set.

	Finally, $\Phi\circ \Psi$ is a continuous self-map of $\mathcal K$, which is convex,
	and its image is contained in a compact convex set.
	Hence by Schauder fixed point theorem, it has a fixed point $u$. Denote $v=\Psi(u)$,
	then $(u,v)$ is solution of the desired PDE system, as asserted.
	By construction, it satisfies the estimate
	\begin{equation}
		e^{-4u}e^{2v}|\alpha_n|^2\leq\eta\,.
	\end{equation}
\end{proof}

\newpage
\bibliographystyle{alpha}
\bibliography{references}
\end{document}